\documentclass[12pt]{article}
\usepackage{latexsym}
\usepackage{graphicx}
\usepackage{cite}
\usepackage{enumitem}
\usepackage{amssymb,mathrsfs,amsmath}
\usepackage{booktabs}
\usepackage[justification=centering]{caption}
\usepackage{subfigure}
\usepackage{float}
\usepackage{color}
\usepackage{arydshln}

\usepackage{mathtools,amsthm}
\usepackage[normalem]{ulem}
\usepackage{mathabx}
\usepackage{acronym}

\usepackage[margin=1.5in]{geometry}%
\newtheorem{theorem}{Theorem}[part]
\newtheorem{definition}{Definition}[part]
\newtheorem{proposition}{Proposition}[part]

\newtheorem{lemma}{Lemma}[part]
\newtheorem{corollary}{Corollary}[part]
\newtheorem{remark}{Remark}[part]
\newtheorem{example}{Example}[part]

\def\x{\mathbf x}
\def\top{\mathsf{T}}

\begin{document}

\title{Unique solvability of weakly homogeneous generalized variational inequalities}

\author{Xueli Bai\thanks{School of Mathematics, Tianjin University,
Tianjin 300072, P.R. China. Email: baixueli@tju.edu.cn.} \and
Zheng-Hai Huang\thanks{Corresponding Author. School of Mathematics,
Tianjin University, Tianjin 300072, P.R. China. Email:
huangzhenghai@tju.edu.cn. Tel:+86-22-27403615 Fax:+86-22-27403615}
\and Mengmeng Zheng\thanks{School of Mathematics, Tianjin University, Tianjin 300072, P.R. China.
Email: zmm941112@tju.edu.cn.}}
\date{}

\maketitle

\begin{abstract}
  An interesting observation is that most pairs of weakly homogeneous mappings have no strongly monotonic property, which is one of the key conditions to ensure the unique solvability of the generalized variational inequality. This paper focuses on studying the unique solvability of the generalized variational inequality with a pair of weakly homogeneous mappings. By using a weaker condition than the strong monotonicity and some additional conditions, we achieve several results on the unique solvability of the underlying problem. These results are exported by making use of the exceptional family of elements or derived from new obtained Karamardian-type theorems or established under the exceptional regularity condition. They are new even when the problem comes down to its important subclasses studied in recent years.\vspace{3mm}

\noindent {\bf Key words:}\hspace{2mm} Generalized variational inequality, Weakly homogeneous mapping, Exceptional family of elements, Degree theory, Strictly monotone mapping.  \vspace{3mm}

\noindent {\bf Mathematics Subject
Classifications (2010):}\hspace{2mm} 65K10, 90C33 \vspace{3mm}
\end{abstract}

\section{Introduction}\label{Int}

Variational inequalities (VIs) and complementarity problems (CPs) have been widely studied because of their applications in many fields (see \cite{FP03,FP032,HXQ06,HP90} for example). The unique solvability of these problems has always been one of the important issues, which has been extensively studied in the literature (see \cite{CPS92,DT95,GP92,GS07,GSR03,KF96,K76,K762,MK77,MH-11,M74} for example).

In 1988, Noor \cite{Noor88} introduced a class of generalized variational inequalities (GVIs), which contains VIs and CPs as subclasses. The unique solvability of the GVI can be guaranteed under several conditions, where one of the key conditions is the {\it strong monotonicity of the mapping pair} involved (see \cite{PY95} for example). The strong monotonicity of a mapping pair is a generalization of the strong monotonicity of a single mapping. The latter is a classical condition to guarantee the unique solvability of VIs (see \cite{FP03}). We find that the pair of weakly homogeneous mappings generally does not have the property of strong monotonicity.

In recent years, several classes of special VIs and CPs have attracted people's attention, including tensor complementarity problems (TCPs) (see \cite{HQ191,HQ193,HQ192}), polynomial complementarity problems (PCPs) (see \cite{G17}), generalized polynomial complementarity problems (PGCPs) (see \cite{LLH-18}), tensor variational inequalities (TVIs) (see \cite{WHQ18}), polynomial variational inequalities (PVIs) (see \cite{H18}), and generalized polynomial variational inequalities (PGVIs) (see \cite{WHX19}). With the help of structural properties of tensors and properties of polynomials, lots of theoretical results for these problems have been obtained. Recently, Gowda and Sossa \cite{GS19} studied the {\it variational inequality with a weakly homogeneous mapping} (WHVI) over a finite dimensional real Hilbert space, which is a unified model for the above classes of special problems. By making use of the degree theory and properties of the weakly homogeneous mapping, they obtained several profound results on the nonemptiness and compactness of solution sets of WHVIs. Moreover, they also obtained a uniquely solvable result of the {\it complementarity problem with a weakly homogeneous mapping} (WHCP), a subclass of WHVIs. More recently, the nonemptiness and compactness of solution sets of WHVIs was also investigated by \cite{MZH19}.

Inspired by the works mentioned above, in this paper, we investigate the unique solvability of the GVI with a pair of weakly homogeneous mappings (WHGVI) over a finite dimensional real Hilbert space. 
The contribution of this paper is threefold. 
\begin{itemize}
\item First, we introduce a definition of exceptional family of elements for a pair of mappings and establish an alternative theorem for the WHGVI, and by which, we show that the WHGVI has a unique solution under some assumptions, where one of the key conditions is the strict monotonicity which is weaker than the strong monotonicity. An example is constructed to claim the advantage of the achieved result. Incidentally, we also get a new result on the nonemptiness and compactness of solution sets of WHGVIs.
\item Second, after obtaining a result on the nonemptiness and compactness of solution sets of WHGVIs, we extend a Karamardian-type theorem obtained by \cite{GS19} for WHVIs to WHGVIs. Furthermore, several uniquely solvable results of WHGVIs are given by making use of the achieved Karamardian-type results and the strict monotonicity of the involved mapping pair.
\item Third, we derive a result on the nonemptiness and compactness of solution sets of WHGVIs under an exceptional regularity condition and some additional conditions, and by which, we achieve a uniquely solvable result of the WHGVI. An example is given to confirm the achieved result.
\end{itemize}

\noindent Moreover, since the WHGVI contains WHVIs (and more, TCPs, PCPs, PGCPs, WHCPs, TVIs, PVIs, and PGVIs) as its subclasses, we reduce our main results to these subclasses, which give some new observations for these subclasses.

This paper is divided into eight parts. In Sect. \ref{s2}, we briefly recall some basic concepts and conclusions in the VI as well as the degree theory. In particular, we give a definition of exceptional family of elements for a pair of mappings and present an alternative theorem by using the exceptional family of elements. In Sect. \ref{s3-0}, we show that many pairs of weakly homogeneous mappings do not possess the strong monotonic property, which are also illustrated by several examples. In Sect. \ref{s3}, we establish a uniquely solvable result of the WHGVI with the help of the exceptional family of elements for a pair of mappings, and we illustrate that this result is different from the well-known result achieved by Pang and Yao \cite{PY95} by an example. In Sect. \ref{s6}, we establish a Karamardian-type theorem for the WHGVI, and further obtain several uniquely solvable results for WHGVIs. In Sect. \ref{s9}, we investigate the unique solvability of the WHGVI under an exceptional regularity condition and some additional conditions. In Sect. \ref{s4}, we reduce our main results to several subcases of WHGVIs and compare the results with those existing ones for these subcases. In Sect. \ref{s5}, we complete this paper via giving some conclusions.

\section{Preliminary}\label{s2}

Throughout this paper, let $H$ be a finite dimensional real Hilbert space with inner product $\langle \cdot,\cdot\rangle$ and norm $\|\cdot\|$, and $C$ be a closed convex cone in $H$. For any nonempty set $\Omega$ in $H$, int$(\Omega)$, $\partial \Omega$ and $\bar{\Omega}$ denote the interior, boundary and closure of $\Omega$, respectively. In addition, for any continuous mapping $g:H\rightarrow H$ and a nonempty set $K$ in $H$, $C\supseteq g^{-1}(K):=\{x\in H\mid g(x)\in K\}$ means that if $g(x^*)\in K$, then $x^*\in C$.

For any $z\in H$ and a closed convex set $K$ in $H$, $\Pi_{K}(z)$ denotes the orthogonal projection of $z$ onto $K$, which is the unique vector $\bar{z}\in H$ satisfying the inequality $\langle y-\bar{z},\bar{z}-z\rangle\geq0$ for all $y\in K$. Besides, as a mapping, $\Pi_{K}(z)$ is nonexpansive, that is, $\|\Pi_{K}(u)-\Pi_{K}(v)\|\leq\|u-v\|$ holds for any $u,v\in H$. For a projection mapping $\Pi_{K}(\cdot)$, we have the following property:
\begin{equation}\label{projection}
0\in K \;\; \mbox{\rm and}\;\; u\in K^{*}\quad\Longrightarrow\quad
\Pi_{K}(-u)=0,
\end{equation}
where $K^*$ denotes the dual cone of $K$ which is defined by $K^*:=\{u\in H\mid u^Tx\geq0,\forall x\in K\}$. We use $\mathcal{N}_{K}(z)$ to denote the normal cone of $K$ at $z$  which is defined by
\begin{eqnarray*}
\mathcal{N}_{K}(z):=\left\{\begin{array}{ll}
\{u\in H\mid u^T(y-z)\leq0,\forall y\in {K}\},\quad & \mbox{\rm if}\; z\in {K},\\
\emptyset,\quad & \mbox{\rm otherwise},
\end{array}\right.
\end{eqnarray*}
and $K^\infty$ to denote the recession cone of $K$ which is defined by
$$
K^{\infty}:=
\left\{u\in H\mid \exists t_k\rightarrow\infty,\exists x^k\in K \text{~such that~} \lim_{k\rightarrow\infty}\frac{x^k}{t_k}=u\right\}.
$$
Then, with the definition of the recession cone $K^\infty$, we have that the mapping
\begin{equation}\label{qe}
\mathcal{K}(t) = tK + K^\infty,\;\; 0\leq t\leq1
\end{equation}
satisfies the following property:
$$
\mathcal{K}(t)=tK + K^\infty=tK~(t\neq 0)~{\mbox{\rm \ and\ }}~\mathcal{K}(0) = K^\infty,
$$
where the first statement comes from the fact that $K^\infty$ is a cone. In \cite{GS19}, the authors obtained the following result:

\begin{lemma}{\rm(\cite{GS19})}
Let $\mathcal{K}(\cdot)$ be defined as \eqref{qe} and $\theta(\cdot,\cdot):H\times [0,1]\rightarrow H$ be continuous. Then, the mapping $(x,t)\mapsto \Pi_{\mathcal{K}(t)}\theta(x,t)$ is continuous.
\end{lemma}

\subsection{Variational inequalities with weakly homogeneous mappings}\label{vi}

A continuous mapping $f:C\rightarrow H$ is said to be positively homogeneous of degree $\delta$ with $\delta\geq 0$, if $f(\lambda x)=\lambda^\delta f(x)$ holds for any $x\in C$ and $\lambda>0$. Now, we recall the definition of the weakly homogeneous mapping.

\begin{definition}\label{wh}\rm{(\cite{GS19})}
A mapping $f:C\rightarrow H$ is called to be weakly homogeneous of degree $\delta$ if $f=h+g$, where $h:C\rightarrow H$ is positively homogeneous of degree $\delta$ and $g:C\rightarrow H$ is continuous and $g(x)=o(\|x\|^{\delta})$ (that is, $\frac{g(x)}{\|x\|^{\delta}}\rightarrow0$) as $\|x\|\rightarrow\infty$ in $C$.
\end{definition}

Some basic properties of weakly homogeneous mappings are given below.
\begin{proposition}\label{whp}{\rm(\cite{GS19})}
Let $f=h+g$ be a weakly homogeneous mapping of degree $\delta>0$. Then, the following statements hold:
\begin{description}
  \item[{\rm(i)}] $h(0)=0$;
  \item[{\rm(ii)}] $\lim_{\lambda\rightarrow\infty}\frac{g(\lambda x)}{\lambda^{\delta}}=0$ for all $x\in C$;
  \item[{\rm(iii)}] $h(x)=\lim_{\lambda\rightarrow\infty}\frac{f(\lambda x)}{\lambda^{\delta}}$ for all $x\in C$;
  \item[{\rm(iv)}] In the representation $f=h+g$, $h$ and $g$ are unique on $C$.
\end{description}
\end{proposition}

According to item (iii), we use $f^{\infty}$ to represent $h$ in $f$ and call it the {\it leading term}. Besides, for a weakly homogeneous mapping $f$ with degree $\delta>0$, it can be easily seen that
$$
\lim_{k\rightarrow\infty}\frac{f(x^k)}{\|x^k\|^{\delta}} =\lim_{k\rightarrow\infty}\left[h\left(\frac{x^k}{\|x^k\|}\right)+\frac{g(x^k)}{\|x^k\|^{\delta}}\right] =h(\bar{x})=f^{\infty}(\bar{x})
$$
holds for all $\|x^k\|\rightarrow\infty$, where $\bar{x}:=\lim_{k\rightarrow\infty}\frac{x^k}{\|x^k\|}$.

Hereafter, we denote two weakly homogeneous mappings $f$ and $g$ by
\begin{equation}\label{weakhomo}
f(x):=f^{\infty}(x)+\bar{f}(x)+p\quad\text{and}\quad g(x):=g^{\infty}(x)+\bar{g}(x)+q,
\end{equation}
where $f^{\infty}$ and $g^{\infty}$ are two leading terms in $f$ and $g$ with degrees $\delta_1>0$ and $\delta_2>0$, respectively; $p$ and $q$ are two constant items in $f$ and $g$, respectively; and $\bar{f}(x)=f(x)-f^{\infty}(x)-p$ and $\bar{g}(x)=g(x)-g^{\infty}(x)-q$. Obviously, $\bar{f}(x)=o(\|x\|^{\delta_1})$ and $\bar{g}(x)=o(\|x\|^{\delta_2})$ as $\|x\|\rightarrow\infty$.

Given a nonempty closed convex set $K$ in $H$ and two continuous mappings $f,g:H\rightarrow H$. The {\it generalized variational inequality}, denoted by GVI$(f,g,K)$, is to find an $x^*\in H$ such that
\begin{equation}\label{gvi}
g(x^*)\in K,\quad \langle f(x^*),y-g(x^*)\rangle\geq0,\quad \forall y\in K.
\end{equation}
When $f:C\rightarrow H$ and $g:H\rightarrow H$ with $g^{-1}(K)\subseteq C$ are weakly homogeneous mappings with degrees $\delta_1>0$ and $\delta_2>0$, respectively, we call the problem \eqref{gvi} to be a {\it weakly homogeneous generalized variational inequality}, which will be investigated in this paper. In the following, we denote this problem by WHGVI$(f,g,K)$ for notational convenience.
\begin{itemize}
\item When $g(x)=x$, WHGVI$(f,g,K)$ reduces to the WHVI, 
    which is to find an $x^*\in K$ such that
    $$\langle f(x^*),y-x^*\rangle\geq0,\quad\forall y\in K.$$
    We denote it by WHVI$(f,K)$.
\item When $K$ is a cone, WHGVI$(f,g,K)$  is equivalent to a complementarity problem, called the {\it weakly homogeneous generalized complementarity problem}, which is to find an $x^*\in H$ such that
    \begin{eqnarray}\label{gncp-add}
    g(x^*)\in K,\quad f(x^*)\in K^*\quad\mbox{\rm and} \quad \langle f(x^*),g(x^*)\rangle=0.
    \end{eqnarray}
    We denote it by WHGCP$(f,g,K)$.
\item Furthermore, if $g(x)=x$, then WHGCP$(f,g,K)$ reduces to the weakly homogeneous conic complementarity problem, which is denoted by WHCP$(f,K)$
\end{itemize}

\begin{remark}
Actually, the WHGVI is a wide class of problems. Except from the above mentioned VIs and CPs, it also includes many other important problems as its special cases. Thus, by studying the properties of WHGVIs, we can directly obtain many good results about these subclasses (please see Sect. \ref{s4} for details).
\end{remark}

For any GVI$(f,g,K)$, we recall that the {\it natural mapping} (see \cite{FP03} for more details) is defined by
\begin{equation}\label{na-m}
(f,g)^{nat}_{K}(x):=g(x)-\Pi_{K}[g(x)-f(x)].
\end{equation}
With the help of the natural mapping and the same technique in Proposition 1.5.8 given in \cite{FP03}, an equivalent reformulation of GVI$(f, g, K)$ can be easily established.

\begin{lemma}\label{nat}
Let $K$ be a closed convex set in $H$, and $f: C \rightarrow H$ and $g: H \rightarrow H$ be two continuous mappings. Then, $x^{*}\in H$ is a solution of GVI$(f, g, K)$ if and only if $(f, g)^{nat}_{K}(x^{*})=0$.
\end{lemma}

Let SOL$(f,g,K)$ denote the solution set of GVI$(f, g, K)$. Then, by Lemma \ref{nat} it follows that $x^*\in \mbox{SOL}(f,g,K)$ if and only if $(f, g)^{nat}_{K}(x^*)=0$.

\begin{lemma}{\rm (\cite{FP03})}\label{extension}
Let $\Phi:S\subseteq\mathbb{R}^n\rightarrow\mathbb{R}^m$ be a continuous mapping defined on the nonempty closed set $S$. A continuous extension $\bar{\Phi}:\mathbb{R}^n\rightarrow\mathbb{R}^m$ exsits such that $\bar{\Phi}(x)=\Phi(x)$ for all $x\in S$.
\end{lemma}

From Lemma \ref{extension} it can be easily deduced that for any weakly homogeneous mapping $f: C\rightarrow H$, there always exist continuous extension $F$ of $f$ from $C$ to $H$. What is more, if WHGVI$(f,g,K)$ satisfies $g^{-1}(K)\subseteq C$, then we have SOL$(f,g,K)=$SOL$(F,g,K)$.

Next, we give the definitions of three classes of mappings, which reduce to the ones in \cite{PY95} when $D=H=\mathbb{R}^n$.

\begin{definition}\label{func}
Let $K$ be a nonempty closed convex subset of $H$, and $f:D\rightarrow H$ and $g:H\rightarrow H$ be two continuous mappings, where $D$ is a nonempty subset of $H$ with $g^{-1}(K)\subseteq D$. $f$ is said to be
\begin{description}
  \item[\mbox{\rm(i)}] monotone with respect to $g$ on $K$ if
      $$
      g(x),g(y)\in K\quad
      \Longrightarrow\quad [f(x)-f(y)]^T[g(x)-g(y)]\geq0;
      $$
  \item[\mbox{\rm(ii)}] strictly monotone with respect to $g$ on $K$ if
      $$
      [g(x),g(y)\in K,\; \mbox{\rm and}\; x\neq y]\quad
      \Longrightarrow\quad [f(x)-f(y)]^T[g(x)-g(y)]>0;
      $$
  \item[\mbox{\rm(iii)}] strongly monotone with respect to $g$ on $K$ if there exists a scalar $c>0$ such that
      $$
      g(x),g(y)\in K\quad \Longrightarrow\quad
      [f(x)-f(y)]^T[g(x)-g(y)]\geq c\|x-y\|^2.
      $$
\end{description}
When $g$ is the identity mapping, we simple call that $f$ is monotone on $K$, strictly monotone on $K$ and strongly monotone on $K$, respectively.
\end{definition}
\vspace{2mm}

From the above definitions, it can be easily seen that if $f$ is strongly monotone with respect to $g$ on $K$, then $f$ must be strictly monotone with respect to $g$ on $K$. However, the converse is not necessarily true. Besides, we have the following result about strictly monotone mappings, whose proof is very simple, and hence, we omit it here.
\begin{lemma}\label{onesol}
Let $K$ be a closed convex set in $H$, and $f: C \rightarrow H$ and $g:H\rightarrow H$ be two weakly homogeneous mappings defined by \eqref{weakhomo} with $g^{-1}(K)\subseteq C$. Suppose that $f$ is strictly monotone with respect to $g$ on $K$. Then, WHGVI$(f, g, K)$ has no more than one solution.
\end{lemma}

\subsection{Degree theory}\label{dt}

The degree theory has been extensively applied to the investigation of VIs and CPs (see \cite{G93,GP94} for example). In this subsection, we recall some basic notations used in the degree theory (readers can also refer to \cite{FP03,L78,OR70}). Let $\Omega$ be a bounded open set in $H$, $\phi:\bar{\Omega}\rightarrow H$ be a continuous mapping, and $b\in H$ satisfying $b\notin \phi(\partial \Omega)$. Then, the topological degree of $\phi$ over $\Omega$ with respect to $b$ is defined, which is an integer and denoted by deg$(\phi,\Omega,b)$.

In addition, if $x^*\in \Omega$ and $\phi(x)=\phi(x^*)$ has a unique solution $x^*$ in $\bar{\Omega}$, then, let $\Omega'$ be any bounded open set containing $x^*$, deg$(\phi,\Omega',\phi(x^*))$ remains a constant, which is called the index of $\phi$ at $x^*$ and denoted by ind$(\phi,x^*)$. Especially, when the continuous mapping $\varphi:H\rightarrow H$ satisfies $\varphi(0)=0$ if and only if $x=0$, then,
$$
\text{ind}(\varphi,0)=\text{deg}(\varphi,\Omega,\varphi(0))=\text{deg}(\varphi,\Omega,0)
$$
holds for any bounded open set $\Omega$ containing 0.

Furthermore, we review the following conclusions.

\begin{lemma}\label{lem2}{\rm (\cite{OR70})}
Let $\Omega$ be an open bounded set in $H$ and $\phi:\bar{\Omega}\rightarrow H$ be continuous. If $b\in H$ with $b\notin \phi(\partial\Omega)$ and deg$(\phi,\Omega,b)\neq0$, then, $\phi(x)=b$ has a solution in $\Omega$.
\end{lemma}

\begin{lemma}\label{lem1}{\rm (\cite{OR70})}
Let $\Omega$ be an open bounded set in $H$ and $\mathcal{H}(x,t):\bar{\Omega}\times[0,1]\rightarrow H$ be continuous. If $b\in H$ with $b\notin\{\mathcal{H}(x,t):x\in\partial{\Omega},t\in[0,1]\}$, then, deg$(\mathcal{H}(\cdot,t),\Omega,b)$ remains a constant as $t$ varies over $[0,1]$.
\end{lemma}

Lemma \ref{lem1} is also known as the homotopy invariance of degree.
%

\subsection{Exceptional family of elements}\label{ef}

It is well-known that the exceptional family of elements is a powerful tool to study the existence of solutions to CPs (see \cite{IBK97,IO98}) and VIs (see \cite{HHF04,ZH99,ZHQ99}). In the following, referring to \cite{HHF04} and \cite{IBK97}, we present a definition of exceptional family of elements for a pair of mappings over a finite dimensional real Hilbert space.\vspace{2mm}

\begin{definition}\label{efe}
Let $f:D\rightarrow H$ and $g:H\rightarrow H$ be two continuous mappings where $D$ is a nonempty set in $H$, and $K$ be a closed convex set in $H$ with $g^{-1}(K)\subseteq D$. A set of points $\{x^r\}\subset D$ is called an exceptional family of elements for the pair $(f,g)$ with respect to any $\hat{x}\in H$, if
\begin{description}
  \item[\rm (i)] $\|x^r\|\rightarrow\infty$ as $r\rightarrow\infty$;
  \item[\rm (ii)] $g(x^r)\in K$ for any $r>0$;
  \item[\rm (iii)] for any $r>\|\Pi_K(\hat{x})\|$, there exists a real number $\alpha_r>0$ such that
  \begin{equation}\label{nc}
  -[f(x^r)+\alpha_r(g(x^r)-\hat{x})]\in\mathcal{N}_K(g(x^r)).
  \end{equation}
\end{description}
\end{definition}

\begin{remark}\label{efoe}
When $D=H=\mathbb{R}^n$ and $g(x)=x$, Definition \ref{efe} reduces to Definition 2.1 in \cite{HHF04}, in which an exceptional family of elements for the mapping $f$ was defined.
\end{remark}

By employing the degree theory, we can establish an alternative theorem for GVI$(f,g,K)$, which is useful in later analysis.

\begin{theorem}\label{at}
Let $K$ be a nonempty closed convex set in $H$, $f,g:H\rightarrow H$ be two continuous mappings, and $\Omega^{\hat{x}}_r:=\{x\in H\mid \|g(x)\|<r\}$ where $r>\|\Pi_K(\hat{x})\|$ for any given $\hat{x}\in H$. Suppose that
\begin{itemize}
  \item[\mbox{\rm($a$)}] the boundedness of $\|g(x)\|$ implies the boundedness of $\|x\|$; and
  \item[\mbox{\rm($b$)}] deg$(g(\cdot),\Omega^{\hat{x}}_r,\Pi_K(\hat{x}))$ is defined and nonzero.
\end{itemize}
Then, there exists either a solution of GVI$(f,g,K)$ or an exceptional family of elements for the pair $(f,g)$ with respect to any given $\hat{x}\in H$.
\end{theorem}
\begin{proof}
The proof is similar to the one in \cite[Theorem 2.2]{HHF04}. We hereby present it for the integrity of the paper. Suppose that GVI$(f,g,K)$ has no solution. We will show that there exists an exceptional family of elements for the pair $(f,g)$ with respect to any given $\hat{x}\in H$. Let homotopy $\mathcal{H}(\cdot,\cdot):H\times[0,1]\rightarrow H$ be defined by
\begin{equation}\label{homo}
\mathcal{H}(x,t):=g(x)-\Pi_K\{t[g(x)-f(x)]+(1-t)\hat{x}\},
\end{equation}
and let
\begin{equation}\label{set}
S_r:=\{x\in H\mid \|g(x)\|<r\},\;\text{where~} r>0.
\end{equation}
First, we show the following result:

{\bf R1}. For any $r>\|\Pi_K(\hat{x})\|$, there exists $x^r\in \partial S_r$ and $t_r\in[0,1]$ such that $\mathcal{H}(x^r,t_r)=0$.\vspace{2mm}

\noindent To this end, we assume that the result {\bf R1} does not hold and derive a contradiction. Suppose that there exists an $\tilde{r}>\|\Pi_K(\hat{x})\|$ such that
$$
0\notin\{\mathcal{H}(x,t):x\in\partial S_{\tilde{r}},t\in[0,1]\}.
$$
Then, by using item $(a)$, the continuity of $\mathcal{H}$ and Lemma \ref{lem1}, we know that deg$(\mathcal{H}(\cdot,t),S_{\tilde{r}},0)$ remains a constant on $[0,1]$. From \eqref{homo} we have $\mathcal{H}(x,0)=g(x)-\Pi_K(\hat{x})$.
Since $\hat{x}$ is an arbitrary given element in $H$ and $\tilde{r}>\|\Pi_K(\hat{x})\|$, from item $(b)$ we obtain that $\text{deg}(\mathcal{H}(x,0),S_{\tilde{r}},0)\neq0$,
and then
$$
\text{deg}(\mathcal{H}(x,1),S_{\tilde{r}},0)=\text{deg}(\mathcal{H}(x,0),S_{\tilde{r}},0)\neq0.
$$
From Lemma \ref{lem2} and the fact that $\mathcal{H}(x,1)=g(x)-\Pi_K[g(x)-f(x)]$, it immediately follows that $\mathcal{H}(x,1)=0$ has a solution. According to Lemma \ref{nat}, this implies that GVI$(f,g,K)$ has a solution, which is a contradiction. Thus, {\bf R1} holds. Then,
\begin{equation}\label{eq1}
g(x^r)=\Pi_K\{t_r[g(x^r)-f(x^r)]+(1-t_r)\hat{x}\}\in K,
\end{equation}
which indicates that
\begin{equation}\label{eq2}
-\{g(x^r)-[t_r(g(x^r)-f(x^r))+(1-t_r)\hat{x}]\}\in\mathcal{N}_K(g(x^r)).
\end{equation}
On one hand, the fact that GVI$(f,g,K)$ has no solution leads to $\mathcal{H}(x,1)\neq0$, and then, $t_r\neq1$ in \eqref{eq1}. On the other hand, from \eqref{set} and {\bf R1} we obtain that $\|g(x^r)\|=r>\|\Pi_K(\hat{x})\|$, which leads to $\mathcal{H}(x,0)\neq0$, and then, $t_r\neq0$ in \eqref{eq1}. These two aspects together give rise to the fact that $t_r\in(0,1)$ in \eqref{eq2}. Denote $\alpha_r:=(1-t_r)/t_r>0$. From \eqref{eq2} we know that
$$
-[f(x^r)+\alpha_r(g(x^r)-\hat{x})]\in\mathcal{N}_K(g(x^r)).
$$
Besides, for $0<r\leq\|\Pi_K(\hat{x})\|$, let $g(x^r)$ be any point in $K$. Then, we have that $g(x^r)\in K$ and $\|g(x^r)\|\rightarrow\infty$ as $r\rightarrow\infty$. So, based on the continuity of $g$, we obtain that $\|x^r\|\rightarrow\infty$ as $r\rightarrow\infty$.

Therefore, $\{x^r\}$ is an exceptional family of elements for the pair $(f,g)$ with respect to $\hat{x}$.
\end{proof}

\begin{corollary}\label{at1}
Given a nonempty closed convex set $K$ in $H$, and two continuous mappings $g:H\rightarrow H$ and $f: C\rightarrow H$. Let $g^{-1}(K)\subseteq C$ and $\Omega^{\hat{x}}_r:=\{x\in H\mid \|g(x)\|<r\}$ where $r>\|\Pi_K(\hat{x})\|$ for any given $\hat{x}\in H$. Suppose that
\begin{itemize}
  \item[\mbox{\rm($a$)}] the boundedness of $\|g(x)\|$ implies the boundedness of $\|x\|$; and
  \item[\mbox{\rm($b$)}] deg$(g(\cdot),\Omega^{\hat{x}}_r,\Pi_K(\hat{x}))$ is defined and nonzero.
\end{itemize}
Then, there exists either a solution of GVI$(f,g,K)$ or an exceptional family of elements for the pair $(f,g)$ with respect to any given $\hat{x}\in H$.
\end{corollary}
\begin{proof}
Suppose that GVI$(f,g,K)$ has no solution. We will show that there exists an exceptional family of elements for the pair $(f,g)$ with respect to any given $\hat{x}\in H$. Let $F$ be any extension of $f$ to $H$, then, it follows that GVI$(F,g,K)$ has no solution. Thus, following the steps in Theorem \ref{at}, we can get an exceptional family of elements $\{x^r\}\subset H$ for the pair $(F,g)$ with respect to any given $\hat{x}\in H$, which satisfies: $\|x^r\|\rightarrow\infty$ as $r\rightarrow\infty$; $g(x^r)\in K$ for any $r>0$; and for any $r>\|\Pi_K(\hat{x})\|$, there exists a real number $\alpha_r>0$ such that $-[F(x^r)+\alpha_r(g(x^r)-\hat{x})]\in\mathcal{N}_K(g(x^r))$.

Since $g(x^r)\in K$ for any $r>0$, we have $x^r\in g^{-1}(K)\subseteq C$, and then, $F(x^r)=f(x^r)$ for any $r>0$. Thus, the set of points $\{x^r\}\subset H$ also satisfies
$-[f(x^r)+\alpha_r(g(x^r)-\hat{x})]\in\mathcal{N}_K(g(x^r))$ for any $r>\|\Pi_K(\hat{x})\|$,
which shows that $\{x^r\}$ is also an exceptional family of elements for the pair $(f,g)$ with respect to any given $\hat{x}\in H$. This completes the proof.
\end{proof}

\section{Discussions of the strong monotonicity}\label{s3-0}

In this paper, our aim is to investigate the unique solvability of WHGVI$(f,g,K)$, where $f$ and $g$ are two weakly homogeneous mappings defined by \eqref{weakhomo}. To see the need for this research, we first recall a well-known uniquely solvable result of GVI$(f,g,K)$ achieved by Pang and Yao in \cite{PY95}, which is stated as follows.

\begin{theorem}\label{old}
Let $K$ be a nonempty closed convex subset of $\mathbb{R}^n$, and $f,g:\mathbb{R}^n\rightarrow\mathbb{R}^n$ be two continuous functions with $g$ being injective. Suppose there exists a vector $z\in g^{-1}(K)$ and positive scalars $\alpha$ and $L$ such that $\|g(x)-g(z)\|\leq L\|x-z\|$ holds for any $x\in g^{-1}(K)$ with $\|x\|\geq\alpha$. If $f$ is strongly monotone with respect to $g$ on $K$, then there exists a unique vector $\bar{x}\in\mathbb{R}^n$ satisfying
$g(x)=\Pi_{K}[g(x)-f(x)]$.
\end{theorem}


From Lemma \ref{nat} we know that the unique vector $\bar{x}$ in Theorem \ref{old} is actually the unique solution of GVI$(f,g,K)$. To obtain the unique solvability of GVI$(f,g,K)$, Theorem \ref{old} requires that the involved pair of mappings satisfies $\|g(x)-g(z)\|\leq L\|x-z\|$ for any $x\in g^{-1}(K)$ with $\|x\|\geq\alpha$ and possesses the strongly monotonic property. However, it can be seen that these two assumptions may not be true in lots of cases when both $f$ and $g$ are weakly homogeneous mappings. In the following, we only show that for many pairs of weakly homogeneous mappings $f$ and $g$, it is impossible that $f$ is strongly monotone with respect to $g$ on $K$.

\begin{proposition}\label{smne1}
Let $K$ be a nonempty closed convex subset of $H$ and $f,g:C\rightarrow H$ be weakly homogeneous mappings defined by \eqref{weakhomo} with degrees $\delta_1>0$ and $\delta_2>0$, respectively. Suppose that $g^{-1}(K)$ is unbounded. If $\delta_1+\delta_2<2$, then $f$ is not strongly monotone with respect to $g$ on $K$.
\end{proposition}
\begin{proof}
Suppose on the contrary that $f$ is strongly monotone with respect to $g$ on $K$. Then, there exists a scalar $c>0$ such that for any given $g(y)\in K$,
$$
\langle f(x)-f(y),g(x)-g(y)\rangle\geq c\|x-y\|^2
$$
holds for any $g(x)\in K$, that is,
\begin{equation}\label{ssmono1}
\langle f^{\infty}(x)+\bar{f}(x)+p-f(y),g^{\infty}(x)+\bar{g}(x)+q-g(y)\rangle\geq c\|x-y\|^2.
\end{equation}
Obviously the degree of the left-hand side of \eqref{ssmono1} is $\delta_1+\delta_2<2$. Dividing both sides of \eqref{ssmono1} by $\|x-y\|^2$, we obtain that
\begin{equation}\label{ssmono3}
\frac{\langle f^{\infty}(x)+\bar{f}(x)+p-f(y),g^{\infty}(x)+\bar{g}(x)+q-g(y)\rangle}{\|x-y\|^2}\geq c.
\end{equation}
By the unboundedness of $g^{-1}(K)$, there exists an unbounded sequence $\{x^k\}$ such that $g(x^k)\in K$ for any $k$. Let $\|x^k\|\rightarrow\infty$, then, the left-hand side of \eqref{ssmono3} tends to $0$, which is a contradiction!

Thus, $f$ is not strongly monotone with respect to $g$ on $K$.
\end{proof}

\begin{remark}\label{aa1}
Suppose that $g(x)=x$. Then, the condition $\delta_1+\delta_2<2$ in Proposition \ref{smne1} reduces to the degree of $f$ is less than one, that is, if the degree of $f$ is less than one, then $f$ is not strongly monotone on $K$.
\end{remark}

Here, we use an example to illustrate Proposition \ref{smne1}.

\begin{example}\label{exm2-1}
Suppose that $H=\mathbb{R}^2$, $C=\mathbb{R}^2_+$, and $K=\{(s,t)^\top\mid s\geq0, t\geq1\}$. We define two weakly homogeneous mappings from $C$ to $H$ by
$$
f(x)=\left(\begin{array}{c} x_1^{1/2}+2\\ x_2^{1/2}\end{array}\right)\quad \mbox{\rm and}\quad g(x)=\left(\begin{array}{c} x_1^{1/3}\\ x_2^{1/3}+1\end{array}\right).$$
\end{example}

In Example \ref{exm2-1},  $\delta_1+\delta_2=1/2+1/3=5/6<1$. Suppose $f$ is strongly monotone with respect to $g$ on $K$. Since $g(0)\in K$, there exists a positive scalar $c>0$ such that for any $g(x)\in K$,
$$
[f(x)-f(0)]^\top[g(x)-g(0)]=x_1^{5/6}+x_2^{5/6}\geq c\|x\|^2.
$$
Dividing the above inequality both sides by $\|x\|^2$ we obtain that
$$
\frac{x_1^{5/6}+x_2^{5/6}}{\|x\|^2}=\frac{1}{\|x\|^{7/6}}h(\tilde{x})\geq c,
$$
where $\tilde{x}=\frac{x}{\|x\|}$ and $h(x)=x_1^{5/6}+x_2^{5/6}$ is a positive homogeneous function with degree $5/6$. Let $\|x\|\rightarrow\infty$, then, the left-hand side of the above inequality tends to $0$, which is a contradiction! Therefore, $f$ is not strongly monotone with respect to $g$ on $K$.

\begin{proposition}\label{smne2}
Let $K$ be a nonempty closed convex subset of $H$ and $f,g:C\rightarrow H$ be two weakly homogeneous mappings. If there exists some $\hat{x}\in C$ satisfying $g(\hat{x})\in K$ such that
$$\frac{\langle f(x)-f(\hat{x}),g(x)-g(\hat{x})\rangle}{\|x-\hat{x}\|^2}\rightarrow0\quad \mbox{\rm as}\;\; x\rightarrow \hat{x},
$$
then $f$ is not strongly monotone with respect to $g$ on $K$.
\end{proposition}
\begin{proof}
Suppose on the contrary that $f$ is strongly monotone with respect to $g$ on $K$. Then, there exists a positive scalar $c>0$ such that for any $g(x),g(y)\in K$,
$$
\langle f(x)-f(y),g(x)-g(y)\rangle\geq c\|x-y\|^2,
$$
and hence, we have
$$
\frac{\langle f(x)-f(\hat{x}),g(x)-g(\hat{x})\rangle}{\|x-\hat{x}\|^2}\geq c.
$$
Let $x\rightarrow \hat{x}$, then, the left-hand side of the above inequality tends to zero, while the right-hand side is a positive constant, which is a contradiction!

Thus, $f$ is not strongly monotone with respect to $g$ on $K$.
\end{proof}

Now, we present an example to illustrate Proposition \ref{smne2}.

\begin{example}\label{nsm}
Let $H=C=\mathbb{R}^2$ and $K=\mathbb{R}^2_+$. We define two weakly homogeneous mappings from $\mathbb{R}^2$ to $\mathbb{R}^2$ by
$$
f(x)=\left(\begin{array}{c} x_1^3+3\\ x_2^3+6\end{array}\right)\quad \mbox{\rm and}\quad g(x)=\left(\begin{array}{c} x_1^4+\cos x_1+1\\ x_2^4+2\end{array}\right).
$$
\end{example}

In Example \ref{nsm}, since $K=\mathbb{R}^2_+$, we may take $\hat{x}=0$, i.e., $g(\hat{x})\in K$. Suppose that $f$ is strongly monotone with respect to $g$ on $K$. Then, since $g(0)\in K$, there exists a positive scalar $c>0$ such that for any $g(x)\in K$,
$$
\langle f(x)-f(0),g(x)-g(0)\rangle=x_1^7+x_2^7+x_1^3(\cos x_1-1)\geq c\|x\|^2.
$$
Dividing both sides of the above inequality by $\|x\|^2$, we have
$$
\frac{x_1^7+x_2^7+x_1^3(\cos x_1-1)}{\|x\|^2}\geq c.
$$
Let $\|x\|\rightarrow0$, then, the left-hand side of the above inequality tends to zero, while the right-hand side is a positive constant, which is a contradiction! Hence, $f$ is not strongly monotone with respect to $g$ on $K$.\vspace{2mm}

From Proposition \ref{smne2}, the following result holds immediately.
\begin{corollary}\label{coro-add1}
Let $K$ be a nonempty closed convex subset of $H$ and $f,g:C\rightarrow H$ be two weakly homogeneous mappings. Suppose that $f$ and $g$ are finite sums of homogeneous mappings on $C$ of the forms:
$$
f(x)=h_\nu(x)+h_{\nu-1}(x)+\cdots+h_1(x)+h_0(x),
$$
$$
g(x)=\bar{h}_\omega(x)+\bar{h}_{\omega-1}(x)+\cdots+\bar{h}_1(x)+\bar{h}_0(x),
$$
respectively, where $\nu,\omega > 0$ are integers, $h_i(x)$ and $\bar{h}_j(x)$ are positively homogeneous with degrees $\gamma_i$ and $\beta_j$ on $C$, and $\gamma_\nu>\gamma_{\nu-1}>\cdots>\gamma_1>\gamma_0 = 0$, $\beta_\omega>\beta_{\omega-1}>\cdots>\beta_1>\beta_0 = 0$. If $g(0)\in K$ and $\gamma_1+\beta_1>2$, then, $f$ is not strongly monotone with respect to $g$ on $K$.
\end{corollary}

\begin{remark}\label{aa2}
\rm{(i)} Suppose that $g(x)=x$. Then, the conditions $g(0)\in K$ and $\gamma_1+\beta_1>2$ in Corollary \ref{coro-add1} reduce to $0\in K$ and the degree of $h_1(x)$ is no less than one.

\rm{(ii)} Recall that for any positive integers $m$ and $n$ with $m,n\geq2$, ${\cal A}=(a_{i_1i_2\cdots i_m})$, where $a_{i_1i_2\cdots i_m}\in\mathbb{R}$ for $i_j\in\{1,2,\ldots,n\}$ and $j\in\{1,2,\ldots,m\}$, is called an $m$-th order $n$-dimensional tensor. We denote the set of all $m$-th order $n$-dimensional tensor by $\mathbb{R}^{[m,n]}$. For any ${\cal A}=(a_{i_1i_2\cdots i_m})\in \mathbb{R}^{[m,n]}$ and $x=(x_1,\ldots,x_n)^\top\in\mathbb{R}^n$, we have ${\cal A}x^{m-1}\in\mathbb{R}^n$, whose the $i$th component is given by
\begin{equation*}\label{i1}
\left({\cal A}x^{m-1}\right)_i:=\sum^n_{i_2,\cdots,i_m=1}a_{ii_2\cdots i_m}x_{i_2}\cdots x_{i_m},~~\forall i\in \{1,2,\ldots,n\}.
\end{equation*}
In Corollary \ref{coro-add1}, if both weakly homogeneous mappings $f$ and $g$ are polynomials, which are defined by
\begin{equation}\label{poly}
f(x)=\sum^{m-1}_{k=1}{\cal A}^{(k)}x^{m-k}+a \quad\text{and}\quad g(x)=\sum^{l-1}_{p=1}{\cal B}^{(p)}x^{l-p}+b
\end{equation}
where $({\cal A}^{(1)},\ldots,{\cal A}^{(m-1)})\in \mathbb{R}^{[m,n]}\times \cdots\times\mathbb{R}^{[2,n]}$, $({\cal B}^{(1)},\ldots,{\cal B}^{(l-1)})\in \mathbb{R}^{[l,n]}\times \cdots\times\mathbb{R}^{[2,n]}$, $a\in\mathbb{R}^n$, and $b\in\mathbb{R}^n$,
then, Corollary \ref{coro-add1} reduces to Proposition 1 in \cite{WHX19}.
\end{remark}

Just as the strong monotonicity of the mapping plays a role in the study of VIs, the uniform $P$-property of the mapping is one of the important conditions to guarantee that the complementary problem has a unique solution.
At the end of this section, we give some observations on the concept of the uniform $P$-mapping.

\begin{definition}\label{func-uni-p}
The mapping $f:\mathbb{R}^n_+\rightarrow\mathbb{R}^n$ is said to be a unform $P$-mapping with respective to $g:\mathbb{R}^n\rightarrow\mathbb{R}^n$ on $\mathbb{R}^n_+$, if there exists some $\rho>0$ such that
$$
\max_{i\in \{1,2,\ldots,n\}}[f_{i}(x)-f_{i}(y)][g_{i}(x)-g_{i}(y)]\geq\rho\|{x}-{y}\|^2, \quad \forall g(x),g(y)\in\mathbb{R}^n_+.
$$
If $\mathbb{R}^n_+$ is replaced by $\mathbb{R}^n$, we simple call that $f$ is a unform $P$-mapping with respective to $g$.
\end{definition}

Consider a class of generalized complementarity problems, which is GVI$(f,g,K)$ with $H:=\mathbb{R}^n$ and $K:=\mathbb{R}^n_+$.
Similar to the one in \cite{KF96}, one can show that this problem has a unique solution under the assumption that $f$ is a uniform $P$-mapping with respective to $g$ on $\mathbb{R}^n_+$ and some additional conditions.

When $g$ is the identity mapping, the uniform $P$-property of mapping pair $(f,g)$ reduces to the uniform $P$-property of mapping $f$ which is called that $f$ is a uniform $P$-mapping on $\mathbb{R}^n_+$. Such a property is one of the key conditions to ensure the unique solvability of CPs (see \cite{CPS92,FP03,HXQ06} for example).

\begin{remark}\label{aa202}
In a similar way as those in Propositions \ref{smne1} and \ref{smne2}, it is easy to verify that lots of weakly homogeneous mapping pairs $f:\mathbb{R}^n_+\rightarrow\mathbb{R}^n$ and $g:\mathbb{R}^n\rightarrow\mathbb{R}^n$ do not possess the uniform $P$-property described in Definition \ref{func-uni-p}.
\end{remark}

\section{Uniqueness derived by using the exceptionally family of elements}\label{s3}

From Propositions \ref{smne1} and \ref{smne2}, we can see that many pairs of weakly homogeneous mappings do not satisfy the strongly monotonic property. Thus, Theorem \ref{old} cannot be directly applied to the WHGVI in many cases. In the following, we investigate the unique solvability of the WHGVI under the strict monotonicity and some additional assumptions. We also construct an example to compare our result with the famous uniqueness result stated in Theorem \ref{old} in the case of the both involved mappings being weakly homogeneous.

Before showing the main result, we first define
$$
B:=\{x\in H\mid \|x\|=1\}\quad \mbox{\rm and}\quad R:=\{x\in H\mid g^{\infty}(x)\in K^{\infty}\}.
$$
It is easy to see that for a weakly homogeneous mapping $g:H\rightarrow H$ defined by \eqref{weakhomo} with degree $\delta_2>0$, we have that
$$
g(\lambda x)=\lambda^{\delta_2}g^{\infty}(x)+\bar{g}(\lambda x)+q
$$
holds for all $\lambda>0$. Let $\lambda\rightarrow\infty$, we have $\|\lambda x\|\rightarrow\infty$ and $\|g(\lambda x)\|\rightarrow\infty$. Hence, in this case, the boundedness of $\|g(x)\|$ implies the boundedness of $\|x\|$, which means that the condition (a) in Corollary \ref{at1} holds trivially.

\begin{theorem}\label{main}
Given a nonempty closed convex subset $K$ of $H$, and two weakly homogeneous mappings $f:C\rightarrow H$ and $g:H\rightarrow H$ defined by \eqref{weakhomo} with degrees $\delta_1>0$ and $\delta_2>0$, respectively. Let $g^{-1}(K)\subseteq C$ and $\Omega^{\hat{x}}_r:=\{x\in H\mid \|g(x)\|<r\}$ where $r>\|\Pi_K(\hat{x})\|$ for any given $\hat{x}\in H$. Suppose that deg$(g(\cdot),\Omega^{\hat{x}}_r,\Pi_K(\hat{x}))$ is defined and nonzero, and the following conditions hold:
\begin{description}
  \item[\rm(i)] $f$ is strictly monotone with respect to $g$ on $K$; and
  \item[\rm(ii)] $\langle f^{\infty}(x),g^{\infty}(x)\rangle\neq0$ for any $x\in B\bigcap R$.
\end{description}
Then, WHGVI$(f,g,K)$ has a unique solution.
\end{theorem}
\begin{proof}
First, we show that the solution set of WHGVI$(f,g,K)$ is nonempty. Here, we use the proof by contradiction. Suppose on the contrary that WHGVI$(f,g,K)$ has no solution. Then, from Corollary \ref{at1} we know that there exists an exceptional family of elements $\{x^r\}$ for the pair $(f,g)$ with respect to any $\hat{x}\in H$, which satisfies $g(x^r)\in K$ and $\|x^r\|\rightarrow\infty$ as $r\rightarrow\infty$. Let $\hat{x}=0$, then from \eqref{nc} we obtain that for any $r>\|\Pi_K(0)\|$, there exists a scalar $\alpha_r>0$ such that
$$-[f(x^r)+\alpha_rg(x^r)]\in\mathcal{N}_K(g(x^r)).$$
According to the definition of normal cone, we have that for any $r>\|\Pi_K(0)\|$,
\begin{equation}\label{thm-eq1}
\langle y-g(x^r),f(x^r)+\alpha_rg(x^r)\rangle\geq0,\quad \forall y\in K.
\end{equation}
Dividing both sides of \eqref{thm-eq1} by $\|x^r\|^{\delta_1+\delta_2}$, we obtain that for any $r>\|\Pi_K(0)\|$,
\begin{equation}\label{thm-eq2}
\left\langle\frac{y-g(x^r)}{\|x^r\|^{\delta_2}},\frac{f(x^r)}{\|x^r\|^{\delta_1}}\right\rangle +\alpha_r\|x^r\|^{\delta_2-\delta_1}\left\langle\frac{y-g(x^r)}{\|x^r\|^{\delta_2}}, \frac{g(x^r)}{\|x^r\|^{\delta_2}}\right\rangle\geq0,\quad \forall y\in K.
\end{equation}
From condition (i) and $g(x^r)\in K$ we know that for any given $g(\theta)\in K$ and  any $r>\|\Pi_K(0)\|$,
\begin{equation}\label{78}
\langle f(x^r)-f(\theta),g(x^r)-g(\theta)\rangle>0.
\end{equation}
Let $\tilde{x}^r=\frac{x^r}{\|x^r\|}$ and $\tilde{x}^r\rightarrow\tilde{x}$ as $r\rightarrow\infty$, then,
$$
\lim_{r\rightarrow\infty}\frac{f(x^r)}{\|x^r\|^{\delta_1}}=f^{\infty}(\tilde{x})\quad \text{and}\quad\lim_{r\rightarrow\infty}\frac{g(x^r)}{\|x^r\|^{\delta_2}}=g^{\infty}(\tilde{x})\in K^{\infty}.
$$
Obviously, $\tilde{x}\in B\bigcap R$. Dividing both sides of \eqref{78} by $\|x^r\|^{\delta_1+\delta_2}$ and let $r\rightarrow\infty$, we have $\langle f^{\infty}(\tilde{x}),g^{\infty}(\tilde{x})\rangle\geq 0$. This, together with condition (ii), implies that $\langle f^{\infty}(\tilde{x}),g^{\infty}(\tilde{x})\rangle>0$. Thus, for any fixed $y\in K$,
\begin{equation}\label{thm-eq5}
\lim_{r\rightarrow\infty}\left\langle\frac{y-g(x^r)}{\|x^r\|^{\delta_2}}, \frac{f(x^r)}{\|x^r\|^{\delta_1}}\right\rangle=-\langle f^{\infty}(\tilde{x}),g^{\infty}(\tilde{x})\rangle<0.
\end{equation}
Besides, from item (ii) we can also obtain that $g^{\infty}(\tilde{x})\neq0$, which leads to
$$
\lim_{r\rightarrow\infty}\left\langle\frac{y-g(x^r)}{\|x^r\|^{\delta_2}}, \frac{g(x^r)}{\|x^r\|^{\delta_2}}\right\rangle=-\|g^{\infty}(\tilde{x})\|^2<0
$$
for any fixed $y\in K$.
Hence, for all sufficiently large $r$, we have that for any fixed $y\in K$,
$$
\left\langle\frac{y-g(x^r)}{\|x^r\|^{\delta_2}},\frac{g(x^r)}{\|x^r\|^{\delta_2}}\right\rangle<0.
$$
Thus, it follows from (\ref{thm-eq2}) that for all sufficiently large $r$,
$$
\left\langle\frac{y-g(x^r)}{\|x^r\|^{\delta_2}},\frac{f(x^r)}{\|x^r\|^{\delta_1}}\right\rangle \geq0
$$
holds for any fixed $y\in K$, which implies that
$$
\lim_{r\rightarrow\infty}\left\langle\frac{y-g(x^r)}{\|x^r\|^{\delta_2}}, \frac{f(x^r)}{\|x^r\|^{\delta_1}}\right\rangle\geq 0.
$$
This contradicts \eqref{thm-eq5}. Thus, WHGVI$(f,g,K)$ has a nonempty solution set.

From Lemma \ref{onesol}, under the assumption of strict monotonicity, WHGVI$(f,g,K)$ has no more than one solution. Thus, WHGVI$(f,g,K)$ has a unique solution.
\end{proof}

Now, we use other restrictions on the mapping $g$ to replace the degree condition used in Theorem \ref{main}, and get the following result.

\begin{theorem}\label{bi}
Given a nonempty closed convex subset $K$ of $H$, and two weakly homogeneous mappings $f:C\rightarrow H$ and $g:H\rightarrow H$ defined by \eqref{weakhomo} with degrees $\delta_1>0$ and $\delta_2>0$, respectively. Suppose that $g$ is an injective mapping, $g^{-1}(K)\subseteq C$, and for any $y\in K$, there exists an $x\in H$ such that $g(x)=y$. If the following conditions hold:
\begin{description}
  \item[\rm(i)] $f$ is strictly monotone with respect to $g$ on $K$; and
  \item[\rm(ii)] $\langle f^{\infty}(x),g^{\infty}(x)\rangle\neq0$ for any $x\in B\bigcap R$,
\end{description}
then, WHGVI$(f,g,K)$ has a unique solution.
\end{theorem}


Now, we construct an example in which all the conditions in Theorem \ref{bi} are satisfied, but the conditions in Theorem \ref{old} are not satisfied.

\begin{example}\label{ex1}
Let $C=\{(c,0)^{\top}\mid c\geq0\}\subseteq H=\mathbb{R}^2$ and $K=\{(s,0)^{\top}\mid s\geq 2\}$. Consider WHGVI$(f,g,K)$, where $f:C\rightarrow H$ and $g:H\rightarrow H$ are defined as follows:
$$
f(x)=\left(\begin{array}{c}
x_1^{17/3}+x_1^{8/3}x_2^{5/3}+2 \\ x_{2}^{17/3}+x_1^{4/3}x_2+1
\end{array}\right)\quad\mbox{\rm and}\quad
g(x)=\left(\begin{array}{c}
x_{1}^3+\frac{1}{1+x_2^2}+1 \\ x_{2}^3
\end{array}\right).
$$
\end{example}

From Example \ref{ex1}, obviously, $g$ is a continuous injection on $H$ and satisfies $g^{-1}(K)\subseteq C$. Besides, it is also easy to see that
\begin{equation*}
\begin{array}{lll}
f^{\infty}(x)=\left(\begin{array}{c} x_1^{17/3}\\ x_{2}^{17/3}\end{array}\right),&\quad \bar{f}(x)=\left(\begin{array}{c}x_1^{8/3}x_2^{5/3}\\ x_1^{4/3}x_2\end{array}\right),&\quad p=\left(\begin{array}{c}2\\1\end{array}\right),\\
g^{\infty}(x)=\left(\begin{array}{c}x_1^3\\x_2^3\end{array}\right),&\quad \bar{g}(x)=\left(\begin{array}{c}\frac{1}{1+x_2^2}\\0\end{array}\right),&\quad q=\left(\begin{array}{c}1\\0\end{array}\right),
\end{array}
\end{equation*}
and
$$
g^{-1}(K)=\{(x_1,x_2)^\top\mid x_1\geq0,x_2=0\}.
$$

First, for any $g(x),g(y)\in K$, where $x_2=y_2=0$ and $x_1\geq0$, $y_1\geq0$ with $x_1\neq y_1$, we have
\begin{equation*}
[f(x)-f(y)]^\top[g(x)-g(y)]=(x_1^{\frac{17}{3}}-y_1^{\frac{17}{3}})(x^3_1-y^3_1)>0,
\end{equation*}
which means that $f$ is strictly monotone with respect to $g$ on $K$. Besides, for any $x\in B\bigcap R$, we have
$$
\langle f^{\infty}(x),g^{\infty}(x)\rangle=x_1^{\frac{26}{3}}+x_2^{\frac{26}{3}}\neq0.
$$
Thus, all the conditions of Theorem \ref{bi} hold.

Second, we show that the conditions in Theorem \ref{old} are not satisfied. Obviously, from Proposition \ref{smne2}, $f$ is not strongly monotone with respect to $g$ on $K$. Moreover, suppose there exist positive scalars $\alpha>0$ and $L>0$ and a vector $z\in g^{-1}(K)$ such that for all $x\in g^{-1}(K)$ with $\|x\|\geq\alpha$,
$$
\|g(x)-g(z)\|\leq L\|x-z\|.
$$
Since $x,z\in g^{-1}(K)$, we know that $x_2=z_2=0$ and $x_1\geq0$, $z_1\geq 0$. Thus,
$$
\|g(x)-g(z)\|=|x_1^3-z_1^3|\leq L|x_1-z_1|,
$$
which implies that
$$
|x_1^2+x_1z_1+z_1^2|\leq L.
$$
Let $x_1\rightarrow+\infty$, then the left-hand side of the above inequality tends to positive infinity, which is a contradiction! Hence, for WHGVI$(f,g,K)$ in Example \ref{ex1},  the conditions in Theorem \ref{old} are not satisfied.

Last, we show that WHGVI$(f,g,K)$ does have a unique solution. For WHGVI$(f,g,K)$, our purpose is to find $x=(x_1,x_2)^\top$, with $x_1\geq0$ and $x_2=0$, such that
$$
\left(\begin{array}{c}
x_1^{17/3}+x_1^{8/3}x_2^{5/3}+2 \\ x_2^{17/3}+x_1^{4/3}x_2+1
\end{array}\right)^\top
\left(\begin{array}{c}
y_1-x_1^3-\frac{1}{1+x_2^2}-1 \\ y_2-x_2^3
\end{array}\right)\geq0,\quad\forall y_1\geq2,y_2=0,
$$
that is,
$$
\left(x_1^{17/3}+2\right)\left(y_1-x_1^3-2\right)\geq0,\quad\forall y_1\geq 2.
$$
Obviously, $x^*=(0,0)^\top$ is the unique solution of WHGVI$(f,g,K)$ in Example \ref{ex1}.

If the strict monotonicity assumption in Theorem \ref{main} is replaced by a weaker one (i.e., the mapping $f$ is monotone with respect to some fixed vector $g(\theta)\in K$), we can still get the existence of solutions to WHGVI$(f,g,K)$. Furthermore, combining with other conditions in Theorem \ref{main}, we can obtain the compactness of solution sets. This is given as follows.

\begin{theorem}\label{GWHVI}
Given a nonempty closed convex subset $K$ of $H$, and two weakly homogeneous mappings $f:C\rightarrow H$ and $g:H\rightarrow H$ defined by \eqref{weakhomo} with degrees $\delta_1>0$ and $\delta_2>0$, respectively. Let $g^{-1}(K)\subseteq C$ and $\Omega^{\hat{x}}_r:=\{x\in H\mid \|g(x)\|<r\}$ where $r>\|\Pi_K(\hat{x})\|$ for any given $\hat{x}\in H$. Suppose that deg$(g(\cdot),\Omega^{\hat{x}}_r,\Pi_K(\hat{x}))$ is defined and nonzero, and the following conditions hold:
\begin{description}
  \item[\rm (i)] there exists some $g(\theta)\in K$ such that $\langle f(x)-f(\theta),g(x)-g(\theta)\rangle\geq0$ holds for any $g(x)\in K$;
  \item[\rm (ii)] $\langle f^{\infty}(x),g^{\infty}(x)\rangle\neq0$ for any $x\in B\bigcap R$.
\end{description}
Then, WHGVI$(f,g,K)$ has a nonempty compact solution set SOL$(f,g,K)$.
\end{theorem}
\begin{proof}
Following the steps in Theorem \ref{main}, we can easily get that SOL$(f,g,K)$ is nonempty. Now we show the boundedness of SOL$(f,g,K)$. Suppose on the contrary that SOL$(f,g,K)$ is unbounded, then, there exists an unbounded sequence $\{x^k\}\subseteq \text{SOL}(f,g,K)$. Thus, we have
$$
g(x^k)\in K\quad \mbox{\rm and}\quad \langle f(x^k), y-g(x^k)\rangle\geq0,\quad\forall y\in K.
$$
For any $u\in K^{\infty}$ and fixed $g(x^0)\in K$, we have $g(x^0)+\|x^k\|^{\delta_2}u\in K$. By dividing both sides of the above inequality by $\|x^k\|^{\delta_1+\delta_2}$ and taking $y:=g(x^0)+\|x^k\|^{\delta_2}u$, we obtain that
$$
\left\langle\frac{f(x^k)}{\|x^k\|^{\delta_1}},u+\frac{g(x^0)-g(x^k)}{\|x^k\|^{\delta_2}}\right\rangle\geq0,\quad\forall u\in K^{\infty}.
$$
Let $k\rightarrow\infty$ and $\lim_{k\rightarrow\infty}\frac{x^k}{\|x^k\|}=\bar{x}$, then,
\begin{eqnarray}\label{V}
\langle f^{\infty}(\bar{x}),u-g^{\infty}(\bar{x})\rangle\geq0,\quad \forall u\in K^{\infty}.
\end{eqnarray}
According to the definition of the recession cone, we know that
\begin{eqnarray}\label{Vadd1}
g^{\infty}(\bar{x})=\lim_{k\rightarrow\infty}\frac{g(x^k)}{\|x^k\|^{\delta_2}}\in K^{\infty}.
\end{eqnarray}
Since $K^{\infty}$ is a cone, it is easy to obtain from \eqref{V} and \eqref{Vadd1} that $\langle f^{\infty}(\bar{x}),g^{\infty}(\bar{x})\rangle=0$, which is a contradiction to condition (ii)! Therefore, SOL$(f,g,K)$ is bounded. In addition, the closedness of SOL$(f,g,K)$ can be obtained by the continuity of the involved mappings.

Thus, WHGVI$(f,g,K)$ has a nonempty compact solution set.
\end{proof}

By weakening the strict monotonicity assumption in Theorem \ref{main}, we obtained a result on the nonemptines and compactness of the solution set of WHGVI$(f,g,K)$ in Theorem \ref{GWHVI}. However, the conditions of Theorem \ref{GWHVI} cannot guarantee the uniqueness of solutions to WHGVI$(f,g,K)$, which can be seen from the following example.

\begin{example}\label{compact}
Let $H=C=\mathbb{R}^2$ and $K=\{(s,0)^\top\mid s\geq-1\}$. Consider WHGVI$(f,g,K)$, where $f:C\rightarrow H$ and $g:H\rightarrow H$ are defined as follows:
\begin{eqnarray*}
f(x)=\left(
\begin{array}{c}
x_1^3-x_1^2\\
x_2^3+x_2
\end{array}\right)\quad \mbox{\rm and}\quad
g(x)=\left(
\begin{array}{c}
x_1^{3}\\
x_2^{3}
\end{array}\right).
\end{eqnarray*}
\end{example}

In Example \ref{compact}, take $y:=(1,0)^\top$ and $z:=(0,0)^\top$, then
$\langle f(y)-f(z),g(y)-g(z)\rangle=0$,
which indicates that $f$ is not strictly monotone with respect to $g$ on $K$. However, it is easy to see that
\begin{itemize}
\item the degree condition of $g$  in Theorem \ref{GWHVI} holds;
\item take $\theta:=(1,0)^\top$, then $g(\theta)=(1,0)^\top\in K$, and for any $g(x)\in K$, we have $$x_1\geq-1,x_2=0,\quad \mbox{\rm and}\quad \langle f(x)-f(\theta),g(x)-g(\theta)\rangle=x_1^2(x_1-1)^2(x_1^2+x_1+1)\geq0;$$
\item $\langle f^{\infty}(x),g^{\infty}(x)\rangle=x_1^6+x_2^6\neq0$ for any $x\neq0$.
\end{itemize}
Thus, all the conditions in Theorem \ref{GWHVI} are satisfied. By Theorem \ref{GWHVI}, we obtain that the solution set of WHGVI$(f,g,K)$ in Example \ref{compact} is nonempty and compact. However, the solution is not unique. In fact, it is easy to verify that both $x^*=(0,0)^\top$  and $x^*=(1,0)^\top$ are solutions to WHGVI$(f,g,K)$ in Example \ref{compact}.

\section{Uniqueness derived from Karamardian-type theorems}\label{s6}

In \cite{GS19}, the authors established many good theoretical results on the nonemptiness and compactness of solution sets of WHVIs, including a Karamardian-type theorem. First, we generalize one of main results on the nonemptiness and compactness of solution sets of WHVIs in \cite{GS19} to WHGVIs. A uniqueness result is obtained directly.

\begin{theorem}\label{Theorem-2.3}
Given a nonempty closed subset $K$ in $C$, and two weakly homogeneous mappings $f: C\rightarrow H$ and $g : H\rightarrow H$ defined by \eqref{weakhomo} with degrees $\delta_1> 0$ and $\delta_2> 0$, respectively. Let $F$ and $F^\infty$ be any given continuous extensions of $f$ and $f^\infty$, respectively. Suppose that $g$ satisfies that $g^{-1}(C)\subseteq C$ and $\bar{g}(x)+q\in C$ as $\|x\|\rightarrow \infty$, and the following conditions hold:
\begin{description}
\item[\rm ($a$)] SOL$(f^{\infty},g^{\infty},K^\infty)=\{0\}$; and
\item[\rm ($b$)] {\rm ind}$\left((F^{\infty},g^{\infty})^{nat}_{K^\infty}(x),0\right)\neq 0$, where $(F^{\infty},g^{\infty})^{nat}_{K^\infty}(\cdot)$ is defined in \eqref{na-m}.
\end{description}
Then, WHGVI$(f,g,K)$ and WHGCP$(f,g,K^\infty)$ have nonempty compact solution sets.
\end{theorem}
\begin{proof}
In the following, we only show the nonemptiness and compactness of solution sets of WHGVI$(f,g,K)$, since the nonemptiness and compactness of solution sets of WHGCP$(f,g,K^\infty)$ can be obtained by similar steps.

Consider the following homotopy mapping:
\begin{equation*}
\begin{split}
\mathcal{H}(x,t):=&[(1-t)g^\infty(x)+tg(x)]-\\
&~~~~\Pi_{\mathcal{K}(t)}\{[(1-t)g^\infty(x)+tg(x)]-
[(1-t)F^\infty(x)+tF(x)]\},
\end{split}
\end{equation*}
where $\mathcal{K}(t)$ is defined in \eqref{qe}. Then,
$$
\mathcal{H}(x,1)=g(x)-\Pi_K[g(x)-F(x)]\;\;\;{\mbox{\rm and}}\;\;\;\mathcal{H}(x,0)=g^\infty(x)-\Pi_{K^\infty}[g^\infty(x)-F^\infty(x)].
$$
Denote the set of zeros of $\mathcal{H}(x,t)$ by:
$$
\mathbb{Z}:=\{x\in H\mid\mathcal{H}(x,t)=0\; \text{for some}\ t\in [0,1]\}.
$$

Next, we show that $\mathbb{Z}$ is uniformly bounded. For the sake of contradiction, assume that $\mathbb{Z}$ is not uniformly bounded. Then we can find sequences $\{t_{k}\}\subseteq [0,1]$ and $\{0\neq x^{k}\}\subseteq H$
such that\ $\mathcal{H}(x^{k},t_{k})=0$ for any $k$ and $\|x^{k}\|\rightarrow \infty$. It follows from $\mathcal{H}(x^{k},t_{k})=0$ that
\begin{equation*}
\begin{split}
&(1-t_k)g^\infty(x^k)+t_kg(x^k)=\\
&~~~~~\Pi_{\mathcal{K}(t_k)}
\left([(1-t_k)g^\infty(x^k)+t_kg(x^k)]-[(1-t_k)F^\infty(x^k)+t_kF(x^k)]\right),
\end{split}
\end{equation*}
which means that $(1-t_k)g^\infty(x^k)+t_kg(x^k)\in \mathcal{K}(t_k)\subseteq C$ and
\begin{eqnarray}\label{*99}
\left\langle(1-t_k)F^\infty(x^k)+t_kF(x^k),z-[(1-t_k)g^\infty(x^k)+t_kg(x^k)]\right\rangle\geq 0
\end{eqnarray}
holds for any $z\in \mathcal{K}(t_k)$. Since $\bar{g}(x)+q\in C$ as $\|x\|\rightarrow \infty$ and $C$ is a convex cone, it follows that $(1-t_k)[\bar{g}(x^k)+q]\in C$ for all $k$. Furthermore, we have that for any $k$,
$$g(x^k)=(1-t_k)[\bar{g}(x^k)+q]+(1-t_k)g^\infty(x^k)+t_kg(x^k)\in C,$$
which, together with $g^{-1}(C)\subseteq C$, implies that $x^k\in C$. Thus we have that $F(x^k)=f(x^k)$ and $F^\infty(x^k)=f^\infty(x^k)$ for all $k$. Thereby, \eqref{*99} can be written as
$$
\left\langle(1-t_k)f^\infty(x^k)+t_kf(x^k),z-[(1-t_k)g^\infty(x^k)+t_kg(x^k)]\right\rangle\geq 0, \quad \forall \ z\in \mathcal{K}(t_k).
$$
Noting that for any $u\in K^\infty$ and (fixed) $g(x^0)\in K$,  $t_kg(x^0)+\|x^{k}\|^{\delta_2}u\in \mathcal{K}(t_k)$ holds for all $k$. Then, by
choosing $z=t_kg(x^0)+\|x^{k}\|^{\delta_2}u$ and dividing the above relation by $\|x^{k}\|^{\delta_1+\delta_2}$, we get that
\begin{eqnarray}\label{4.1}
\left\langle\frac{(1-t_k)f^\infty(x^k)+t_kf(x^k)}
{\|x^{k}\|^{\delta_1}}, u-\frac{(1-t_k)g^\infty(x^k)+t_kg(x^k)-t_kg(x^0)}{\|x^{k}\|^{\delta_2}}\right\rangle\geq 0\;\;\;
\end{eqnarray}
holds for any $u\in K^\infty$ and all $k$.
Since $\{t_{k}\}$ and $\{\frac{x^{k}}{\|x^{k}\|}\}$ are bounded, we can assume that
\begin{eqnarray*}\label{lim-tx}
\lim_{k\rightarrow\infty}t_{k}=\overline{t} \quad\mbox{\rm and}\quad
\lim_{k\rightarrow\infty}\frac{x^{k}}{\|x^{k}\|}=\overline{x}.
\end{eqnarray*}
Then, let $k\rightarrow\infty$ in \eqref{4.1}, we have that
\begin{eqnarray}\label{4.2}
\langle f^\infty(\overline{x}),u-
g^\infty(\overline{x})\rangle\geq 0, \ \forall\ u\in K^\infty.
\end{eqnarray}
{In the following, we show that $g^\infty(\overline{x})\in K^\infty$.} Noting that $(1-t_k)g^\infty(x^k)+t_kg(x^k)\in \mathcal{K}(t_k)=t_kK+K^\infty$ for all $k$, if $t_k=0$ for infinitely many $k$, then we can find a subsequence $\{x^{k'}\}$ of $\{x^{k}\}$ such that $t_{k'}=0$, thus we can easily get that $$g^\infty(\overline{x}) =\lim_{k'\rightarrow\infty}\frac{(1-t_{k'})g^\infty(x^{k'})+t_{k'}g(x^{k'})}{\|x^{k'}\|^{\delta_2}}\in K^\infty.$$
Otherwise, there must exist infinitely many $k$ such that $t_k>0$, then we can find a subsequence $\{x^{k'}\}$ of $\{x^{k}\}$ such that $t_{k'}>0$. Since $(1-t_{k'})g^\infty(x^{k'})+t_{k'}g(x^{k'})\in \mathcal{K}(t_{k'})=t_{k'}K+K^\infty=t_{k'}K$, each $(1-t_{k'})g^\infty(x^{k'})+t_{k'}g(x^{k'})$ can be written as $(1-t_{k'})g^\infty(x^{k'})+t_{k'}g(x^{k'})=t_{k'}y^{k'}$ with $y^{k'}\in K$.
Noting that $y^{k'}\in K$ and $\frac{\|x^{k'}\|}{t_{k'}}\rightarrow\infty$ as $k'\rightarrow\infty$, then by the definition of the recession cone, we have that
\begin{equation*}
\begin{split}
g^\infty(\overline{x})&=\lim_{k'\rightarrow\infty}\frac{(1-t_{k'})g^\infty(x^{k'})+t_{k'}g(x^{k'})}{\|x^{k'}\|^{\delta_2}}\\
&=\lim_{k'\rightarrow\infty}\frac{t_{k'}y^{k'}}{\|x^{k'}\|^{\delta_2}}=\lim_{k'\rightarrow\infty}\frac{y^{k'}}{\|x^{k'}\|^{\delta_2}/t_{k'}}\in K^\infty.
\end{split}
\end{equation*}
{Therefore, both cases imply that $g^\infty(\overline{x})\in K^\infty$.} This, together with \eqref{4.2}, implies that $0\neq\overline{x}\in \mbox{SOL}(f^{\infty},g^{\infty},K^\infty)$. This is a contradiction! Thereby, $\mathbb{Z}$ is uniformly bounded.

Now, let $\Omega$ be a bounded open set in $H$, which contains $\mathbb{Z}$, then, $0\notin \mathcal{H}(\partial \Omega,t)$ for any $t\in [0,1]$. By the homotopy invariance principle of the degree, we have that,
$$
\text{deg}(\mathcal{H}(\mathbf{x},1),\Omega,0)=\text{deg}(\mathcal{H}(x,0),
\Omega,0)=\text{ind}\left((F^{\infty},g^{\infty})^{nat}_{K^\infty}(x),0\right)\neq 0.
$$
Hence, from Lemma \ref{lem2}, we obtain that $\mbox{SOL}(f,g,K)$ is nonempty. In addition, it follows that SOL$(f,g, K)$ is bounded from the boundedness of $\mathbb{Z}$ which contains SOL$(f,g,K)$. Moreover, the closeness of SOL$(f,g,K)$ is obvious. Therefore, SOL$(f,g,K)$ is nonempty and compact.
\end{proof}

\begin{remark}
When $g(x)=g^\infty(x)=x$, WHGVI$(f,g,K)$ and WHGCP$(f,g,K^\infty)$ reduce to WHVI$(f,K)$ and WHCP$(f,K^\infty)$, respectively. In this case, it is obvious that $g^{-1}(C)=C$ and $\bar{g}(x)+q=0\in C$ for any $x$. Thus, Theorem \ref{Theorem-2.3} can reduce to  \cite[Theorem 4.1]{GS19}.
\end{remark}

\begin{corollary}\label{c33333}
Suppose that all the conditions in Theorem \ref{Theorem-2.3} are satisfied.
If additionally $f$ is strictly monotone with respect to $g$ on $K(K^\infty)$, then, WHGVI$(f,g,K)$ (WHGCP$(f,g,K^\infty)$) has a unique solution.
\end{corollary}

Second, we establish a Karamardian-type theorem for the WHGVI, which is a generalization of the one in \cite{GS19}. Two uniqueness results are also given.
\begin{theorem}\label{Theorem 2.4}
Given a nonempty closed subset $K$ of $C$ with $K^\infty$ being pointed, and two weakly homogeneous mappings $f: C\rightarrow H$ and $g : H\rightarrow H$ defined by \eqref{weakhomo} with degrees $\delta_1> 0$ and $\delta_2> 0$, respectively. Let $F$ and $F^\infty$ be any given continuous extensions of $f^\infty+\bar{f}-\bar{f}(0)$ and $f^\infty$, respectively. Suppose that $g$ satisfies $g^{-1}(C)\subseteq C$ and $\bar{g}(x)+q\in C$ as $\|x\|\rightarrow \infty$, and there exists a vector $d\in \mbox{int}((K^\infty)^*)$ such that one of the following conditions holds:
\begin{description}
\item[\rm ($a$)] SOL$(f^{\infty},g^{\infty},K^\infty)=\{0\}=$SOL$(f^\infty+\bar{f}-\bar{f}(0)+d,g,K^\infty)$ and there exists a nonempty bounded open set $\Omega$ satisfying
    \begin{equation}\label{dist1}
    \text{dist}(F(x),g(x))<\text{dist}(d, \partial((K^\infty)^{*})),\quad \forall x\in \Omega
    \end{equation}
    such that deg$(g(x), \Omega, 0)\neq 0$.
\item[\rm ($b$)] SOL$(f^{\infty},g^{\infty},K^\infty)=\{0\}=$SOL$(f^\infty+d,g^\infty,K^\infty)$ and there exists a nonempty bounded open set $\Omega$ satisfying
    \begin{equation}\label{dist2}
    \text{dist}(F^\infty(x),g^\infty(x))<\text{dist}(d, \partial((K^\infty)^{*})),\quad \forall x\in \Omega
    \end{equation}
    such that deg$(g^\infty(x), \Omega, 0)\neq 0$.
\end{description}
Then, WHGVI$(f,g,K)$ and WHGCP$(f,g,K^\infty)$ have nonempty compact solution sets.
\end{theorem}
\begin{proof}
By Theorem \ref{Theorem-2.3}, we only need to show ind$\left((F^{\infty},g^{\infty})^{nat}_{K^\infty}(x),0\right)\neq 0$ under the condition $(a)$ or $(b)$.

{\bf Case 1}: Suppose that $(a)$ holds. We consider the following homotopy mapping:
\begin{equation*}
\begin{split}
\mathcal{H}(x,t):=&[(1-t)g^\infty(x)+tg(x)]-\\
&~~~\Pi_{K^\infty}\{[(1-t)g^\infty(x)+tg(x)]-[(1-t)F^\infty(x)+tF(x)+td]\}.
\end{split}
\end{equation*}
Then,
$$
\mathcal{H}(x,1)=g(x)-\Pi_{K^\infty}[g(x)-(F(x)+d)]\;\text{and}\;\mathcal{H}(x,0)=g^\infty(x)-\Pi_{K^\infty}[g^\infty(x)-F^\infty(x)].
$$
Denote the set of zeros of $\mathcal{H}(x,t)$ by:
$$
\mathbb{Z}:=\{x\in H\mid\mathcal{H}(x,t)=0\;\; \mbox{\rm for some}\ t\in [0,1]\}.
$$

By a similar method to the proof of Theorem \ref{Theorem-2.3}, we can obtain that $\mathbb{Z}$ is uniformly bounded, thus SOL$(f,g,K)$ is bounded. Let $\Omega'$ be a bounded open set in $H$, which contains $\mathbb{Z}$, then, $0\notin \mathcal{H}(\partial\Omega',t)$ for any $t\in [0,1]$. Furthermore, by the homotopy invariance principle of the degree, we have that,
\begin{eqnarray}\label{4.4}
\text{ind}\left((F^{\infty},g^{\infty})^{nat}_{K^\infty}(x),0\right)=\text{deg}(\mathcal{H}(x,0),
\Omega',0)=\text{deg}(\mathcal{H}(x,1),\Omega',0).
\end{eqnarray}
It follows from condition $(a)$ that $[F(x)+d]-g(x)$ is concluded in some neighbourhood of $d$ in $(K^\infty)^{*}$ for any $x\in \Omega$. Thus for any $x\in \Omega$,
$$
\mathcal{H}(x,1)=g(x)-\Pi_{K^\infty}\left(-[(F(x)+d)-g(x)]\right)=g(x)-0=g(x).
$$
Thus deg$(\mathcal{H}(x,1),\Omega,0)=$deg$(g(x), \Omega, 0)\neq 0$. Furthermore, by condition $(a)$,
$$
\text{deg}(\mathcal{H}(x,1),\Omega',0)=\text{ind}(\mathcal{H}(x,1),0) =\text{deg}(\mathcal{H}(x,1),\Omega,0)\neq 0.
$$
Therefore, by \eqref{4.4}, we get that ind$\left((F^{\infty},g^{\infty})^{nat}_{K^\infty}(x),0\right)\neq 0$. Hence, from Lemma \ref{lem2}, we obtain that SOL$(f,g,K)$ is nonempty. Above, we have shown that SOL$(f,g,K)$ is nonempty and bounded. It is easy to see that SOL$(f,g,K)$ is closed. Therefore, SOL$(f,g,K)$ is nonempty and compact.

{\bf Case 2}: Suppose that $(b)$ holds. By a similar technique as the one in Case 1, it is not difficult to obtain that SOL$(f,g,K)$ is nonempty and compact.

Thus, either $(a)$ or $(b)$ implies that WHGVI$(f,g,K)$ has a nonempty compact solution set.
\end{proof}

\begin{remark}
When $g(x)=g^\infty(x)=x$, WHGVI$(f,g,K)$ and WHGCP$(f,g,K^\infty)$ reduce to WHVI$(f,K)$ and WHCP$(f,K^\infty)$, respectively. In this case, it is obvious that $g^{-1}(C)=C$, $\bar{g}(x)+q=0\in C$ for any $x$, and deg$(g(x), \Omega, 0)=1=$deg$(g^\infty(x), \Omega, 0)$ for any small open neighbourhood of $0$. Thus, it is not difficult to see that Theorem \ref{Theorem 2.4} can reduce to Theorem 5.1 in \cite{GS19} when WHGVI$(f,g,K)$ reduces to WHVI$(f,K)$.
\end{remark}

\begin{corollary}\label{cor-add1}
Suppose that all the conditions in Theorem \ref{Theorem 2.4} are satisfied.
If additionally $f$ is strictly monotone with respect to $g$ on $K(K^\infty)$, then, WHGVI$(f,g,K)$ (WHGCP$(f,g,K^\infty)$) has a unique solution.
\end{corollary}

In the same way as in Theorem 2.3 in \cite{M74-1}, we can obtain the following result for WHGCPs.
\begin{lemma}\label{onesol2}
Let $f: \mathbb{R}^n_+ \rightarrow \mathbb{R}^n$ and $g:\mathbb{R}^n\rightarrow \mathbb{R}^n$ be two weakly homogeneous mappings defined by \eqref{weakhomo} with $g^{-1}(K)\subseteq \mathbb{R}^n_+$. Suppose that $f$ has $P$-property with respect to $g$ on $\mathbb{R}^n_+$, that is,
$$
\max_{i\in [n]}[f_i(x)-f_i(y)][g_i(x)-g_i(y)]>0, \ \forall\ g(x),g(y)\in \mathbb{R}^n_+\ \mbox{\rm and}\ x\neq y.
$$
Then, WHGCP$(f, g, C)$ has no more than one solution.
\end{lemma}

\begin{remark}\label{cor-add2}
Suppose that $H=\mathbb{R}^n$ and $K=C=\mathbb{R}^n_+$. If the condition that $f$ is strictly monotone with respect to $g$ on $K(K^\infty)$ in Corollaries \ref{c33333} and \ref{cor-add1} is replaced by the condition that $f$ has $P$-property with respect to $g$ on $K$, then the corresponding WHGCP has a unique solution.
\end{remark}

As is known to us, the condition which is described by topological degree is difficult to check directly, in general. Thus, we give the following result where the degree-theoretical conditions are replaced by other properties of mappings.

\begin{theorem}\label{Theorem 2.5}
Given a nonempty closed subset $K$ of $C$ with $K^\infty$ being pointed, and two weakly homogeneous mappings $f: C\rightarrow H$ and $g : H\rightarrow H$ defined by \eqref{weakhomo} with degrees $\delta_1> 0$ and $\delta_2> 0$, respectively. Let $F$ and $F^\infty$ be any given continuous extensions of $f^\infty+\bar{f}-\bar{f}(0)$ and $f^\infty$, respectively.
{Suppose that $g$ satisfies that $g^{-1}(C)\subseteq C$ and $\bar{g}(x)+q\in C$ as $\|x\|\rightarrow \infty$}, and there exists a vector $d\in \mbox{int}((K^\infty)^*)$ such that one of the following conditions holds:
\begin{description}
\item[\rm ($a$)] SOL$(f^{\infty},g^{\infty},K^\infty)=\{0\}=$SOL$(f^\infty+\bar{f}-\bar{f}(0)+d,g,K^\infty)$ and $g$ is an injective mapping satisfying there exists an $x^*$ such that $g(x^*)=0$ and
    $$
    \text{dist}(F(x),g(x))<\text{dist}(d, \partial((K^\infty)^{*}))\;\; \text{as}\;\; x\rightarrow x^*;
    $$
\item[\rm ($b$)] SOL$(f^{\infty},g^{\infty},K^\infty)=\{0\}=$SOL$(f^\infty+d,g^\infty,K^\infty)$ and $g^\infty$ is an injective mapping.
\end{description}
Then, WHGVI$(f,g,K)$ and WHGCP$(f,g,K^\infty)$ have nonempty compact solution sets. Furthermore, if additionally $f$ is strictly monotone with respect to $g$ on $K(K^\infty)$, then WHGVI$(f,g,K)$ (WHGCP$(f,g,K^\infty)$) has a unique solution.
\end{theorem}

\begin{proof}
First, we show that condition $(a)$ can imply condition $(a)$ of Theorem \ref{Theorem 2.4}. If there exists some $x^*$  such that $g(x^*)=0$ and dist$(F(x),g(x))<$dist$(d, \partial((K^\infty)^{*}))$ as $x\rightarrow x^*$, then choosing $\Omega$ being an open neighborhood of $x^*$ such that \eqref{dist1} holds for any $x\in \Omega$, we can obtain that deg$(g(x), \Omega, 0)\neq 0$ by the condition that $g$ is injective.

Second, we show that condition $(b)$ can imply condition $(b)$ of Theorem \ref{Theorem 2.4}. Since $f^\infty(0)=g^\infty(0)=0$, we have that  dist$(F^\infty(x),g^\infty(x))\rightarrow 0$ as $x\rightarrow 0$. Thus choosing $\Omega$ being an small enough open neighborhood of $0$ such that \eqref{dist2} holds for any $x\in \Omega$, we can obtain that deg$(g^\infty(x), \Omega, 0)\neq 0$ by the condition that $g^\infty$ is an injective mapping.
\end{proof}

\begin{remark}
\rm{(i)} Suppose that $g(0)=0$, then the conditions described by the distance inequalities in $(a)$ of Theorems \ref{Theorem 2.4} and \ref{Theorem 2.5} hold naturally.

\rm{(ii)} If $f$ is strictly monotone with respective to $g$ on $H$, then the condition that $g$ is an injective mapping in $(a)$ of Theorem \ref{Theorem 2.5} holds naturally.
\end{remark}

It can be easily seen that the uniqueness result in Theorem \ref{Theorem 2.5} was obtained from a Karamardian-type theorem, and the uniqueness result in Theorem \ref{bi} was derived by using the exceptionally family of elements. Next, we show that the conditions of these two theorems cannot be contained each other, which can be seen from the following examples.

\begin{example}\label{ex-5}
Let $H=\mathbb{R}^3$, $C=\{(x_1,x_2,x_3)^\top\in \mathbb{R}^3\mid x_3\geq0\}$, and $K=\{(x_1,x_2,x_3)^\top\\
\in \mathbb{R}^3\mid x_1=x_2\geq 3,x_3=\pi^3\}$. Consider WHGVI$(f,g,K)$, where $f:C\rightarrow H$ and $g:H\rightarrow H$ are defined as follows:
\begin{eqnarray*}
f(x)=\left(
\begin{array}{c}
x_1^{3}+x_1\sin x_3+1\\
x_2^{3}+x_2\sin x_3+1\\
x_3^3
\end{array}\right)\quad \mbox{\rm and}\quad
g(x)=\left(
\begin{array}{c}
x_1^{3}+\sin(-\frac{\pi}{2}\cdot x_1)+1\\
x_2^{3}+\sin(-\frac{\pi}{2}\cdot x_2)+1\\
x_3^3
\end{array}\right).
\end{eqnarray*}
\end{example}

It is easy to check that $f$ and $g$ are weakly homogeneous mappings with degrees $3$,  $g^{-1}(K)\subseteq C$, $K^\infty=\{(x_1,x_2,x_3)^\top\in \mathbb{R}^3\mid x_1=x_2\geq 0,x_3=0\}$, and $(K^\infty)^*=\{(x_1,x_2,x_3)^\top\in \mathbb{R}^3\mid x_1+x_2\geq 0\}$. Below, we show that for WHGVI$(f,g,K)$ in Example \ref{ex-5}, all of the conditions in Theorem \ref{Theorem 2.5} hold, but at least one of the conditions in Theorem \ref{bi} is not satisfied. We also show that this WHGVI has a unique solution, which conforms the result of Theorem \ref{Theorem 2.5}.

{\it Part I.} We show that all the conditions in Theorem \ref{Theorem 2.5} hold for WHGVI$(f,g,K)$  in Example \ref{ex-5}.
\begin{itemize}
  \item It is not difficult to obtain that $g^{-1}(C)\subseteq C$ and $\bar{g}(x)+q\in C$ as $\|x\|\rightarrow\infty$ by $[\bar{g}(x)+q]_3=0$ for any $x\in \mathbb{R}^3$.
  \item Obviously, $g^\infty$ is an injective mapping on $\mathbb{R}^3$. Next, we discuss the uniqueness of solutions to the corresponding WHGCP$(f^\infty,g^\infty,K^\infty)$ and WHGCP$(f^\infty+d,g^\infty,K^\infty)$ with $d\in \text{int}((K^\infty)^*)$, respectively. Suppose $\bar{x}\in \text{SOL}(f^{\infty},g^{\infty},K^\infty)$, it follows from $g^\infty(\bar{x})\in K^\infty$ and $f^\infty(\bar{x})\in (K^\infty)^*$ that $\bar{x}_1=\bar{x}_2\geq 0$, which, together with $\langle f^\infty(\bar{x}), g^\infty(\bar{x})\rangle=\bar{x}_1^{6}+\bar{x}_2^{6}=0$, implies that $\bar{x}_1=\bar{x}_2=0$, i.e., SOL$(f^{\infty},g^{\infty},K^\infty)=\{0\}$. Similarly, if $\bar{x}\in \text{SOL}(f^{\infty}+d,g^{\infty},K^\infty)$, then it follows from $g^\infty(\bar{x})\in K^\infty$ and $f^\infty(\bar{x})+d\in (K^\infty)^*$ that $\bar{x}_2=\bar{x}_1\geq0$, which, together with
      \begin{equation*}
      \begin{split}
      \langle f^\infty(\bar{x})+d, g^\infty(\bar{x})\rangle&=\bar{x}_1^{6}+d_1\bar{x}_1^3+\bar{x}_2^{6}+d_2\bar{x}_2^3   =2\bar{x}_1^{6}+(d_1+d_2)\bar{x}_1^3=0
      \end{split}
      \end{equation*}
      and $d_1+d_2>0$ (from $d\in \text{int}((K^\infty)^*)$) shows that $\bar{x}_1=\bar{x}_2=0$, i.e., SOL$(f^{\infty}+d,g^{\infty},K^\infty)=\{0\}$. Thereby, the condition $(b)$ in Theorem \ref{Theorem 2.5} holds.
  \item For any $g(x),g(y)\in K$ and $x\neq y$, we have $x_3=y_3=\pi$, $x_1\geq 1, y_1\geq 1$, $x_2\geq 1, y_2\geq 1$ and there exists an index $i\in \{1,2\}$ such that $x_i\neq y_i$. Now, we discuss the strict monotonicity of the mapping $h(x)=x^3+\sin(-\frac{\pi}{2}\cdot x)$ where $1\leq x\in \mathbb{R}$. Since
      $$
      [x^3+\sin(-\frac{\pi}{2}\cdot x)]'=3x^2-\frac{\pi}{2}\cos(-\frac{\pi}{2}\cdot x)>3-\frac{\pi}{2}>0,
      $$
      when $x\geq1$, $h(x)>h(y)$ if $x>y\geq 1$. Thus, for any $x\geq1,y\geq 1$, and $x\neq y$, we can obtain that
      $$
      (x-y)[h(x)-h(y)]=(x-y)[x^3+\sin(-\frac{\pi}{2}\cdot x)-y^3-\sin(-\frac{\pi}{2}\cdot y)]>0.
      $$
      Furthermore, we have that for any $g(x),g(y)\in K$ and $x\neq y$,
      $$
      \begin{array}{rl}
      &[f(x)-f(y)]^\top[g(x)-g(y)]=\\
      &~~~~~\sum_{i=1}^2(x_i^{3}-y_i^{3})[x_i^3+\sin(-\frac{\pi}{2}\cdot x_i)-y_i^3-\sin(-\frac{\pi}{2}\cdot y_i)]>0.
      \end{array}
      $$
      Thereby, $f$ is strictly monotone with respect to $g$ on $K$.
\end{itemize}
Combining the above three cases, we obtain that all the conditions in Theorem \ref{Theorem 2.5} hold for WHGVI$(f,g,K)$ in Example \ref{ex-5}.

{\it Part II.} Since $g(0)=(1,1,0)^\top=g(z)$ where $z=(1,1,0)^\top$, $g$ is not an injective mapping on $H=\mathbb{R}^3$, i.e., at least one of the conditions in Theorem \ref{bi} is not satisfied for WHGVI$(f,g,K)$ in Example \ref{ex-5}.

{\it Part III.} We are going to show that WHGVI$(f,g,K)$ has a unique solution. Considering WHGVI$(f,g,K)$ is to find $x=(x_1,x_2,x_3)^\top$ such that
\begin{eqnarray}\label{G}
g(x)\in K\;\;\mbox{\rm and}\;\;
\left(\begin{array}{c}
x_1^{3}+x_1\sin x_3+1\\
x_2^{3}+x_2\sin x_3+1\\
x_3^3
\end{array}\right)^\top
\left(\begin{array}{c}
y_1-[x_1^{3}+\sin(-\frac{\pi}{2}\cdot x_1)+1]\\
y_2-[x_2^{3}+\sin(-\frac{\pi}{2}\cdot x_2)+1]\\
y_3-x_3^3
\end{array}\right)\geq0\;\;
\end{eqnarray}
for all $y_1\geq3, y_2\geq 3$ and $y_3=\pi^3$. It follows from $g(x)\in K$ that $x_1\geq1$, $x_2\geq1$ and $x_3=\pi$, thus \eqref{G} can be rewritten as
$$
\sum_{i=1}^2(x_i^{3}+1)(y_i-x_i^3-\sin(-\frac{\pi}{2}\cdot x_i)-1)\geq0,\quad\forall y_1\geq3, y_2\geq 3 \;\;\mbox{\rm and}\;\; y_3=\pi^3.
$$
Suppose that $x^*$ is a solution of WHGVI$(f,g,K)$. Then by taking $y=(3,3,\pi^3)^\top$, we have
\begin{eqnarray}\label{5.1}
\sum_{i=1}^2[(x_i^*)^{3}+1][2-(x_i^*)^3-\sin(-\frac{\pi}{2}\cdot x_i^*)]\geq0.
\end{eqnarray}
Since $g(x^*)\in K$, we have that $x_3^*=\pi$, $x_i^*\geq1$ and $(x_i^*)^3+\sin(-\frac{\pi}{2}\cdot x_i^*)\geq 2$ for any $i\in \{1,2\}$, which together with \eqref{5.1} implies that $(x_i^*)^3+\sin(-\frac{\pi}{2}\cdot x_i^*)=2$ for any $i\in \{1,2\}$. Furthermore, noting that $h(x)=x^3+\sin(-\frac{\pi}{2}\cdot x)$ is strictly increasing on $\{x\in \mathbb{R}\mid x\geq 1\}$ and \eqref{G} holds when $x_3^*=\pi$ and $(x_i^*)^3+\sin(-\frac{\pi}{2}\cdot x_i^*)=2$ for any $i\in \{1,2\}$, thus WHGVI$(f,g,K)$ in Example \ref{ex-5} has a unique solution.

At the end of this section, we use Example \ref{ex1} to illustrate the case that all of the conditions in Theorem \ref{bi} are satisfied, but at least one of the conditions in Theorem \ref{Theorem 2.5} is not satisfied. Actually, from the analysis of Example \ref{ex1} given in the above section, we have already known that all the conditions in Theorem \ref{bi} are satisfied for WHGVI$(f,g,K)$ in Example \ref{ex1}. However, since
$$
g^{-1}(C)=\left\{(c',0)^\top\mid c'\geq \sqrt[3]{-2}\right\},
$$
 it is easy to see that $g^{-1}(C)$ is not contained in $C$. Thus, the condition $g^{-1}(C)\subseteq C$ in Theorem \ref{Theorem 2.5} does not hold for WHGVI$(f,g,K)$ in Example \ref{ex1}.

From the above two examples, it can be seen that Theorem \ref{bi} and Theorem \ref{Theorem 2.5} are two different results, which analyze the properties of solutions to WHGVI$(f,g,K)$ from two different aspects.

\section{Uniqueness derived under the exceptional regularity condition}\label{s9}
In Section \ref{s6}, we have obtained a Karamardian-type theorem for the WHGVI, which can reduce to the one established by Gowda and Sossa \cite{GS19} for the WHVI. To the best of our knowledge, the Karamardian-type theorem achieved by \cite{GS19} can cover a lot of existing results obtained recently in TCPs, PCPs, TVIs and PVIs. However, we note that there are some papers which study the properties of solution sets of VIs and CPs by using the exceptional regularity of the involved mappings (see \cite{LLH-18,ZI-00,ZHM-19}). It is not clear whether or not the exceptional regularity of the involved mappings can lead to new results on the nonemptiness and compactness of solution sets and/or the uniqueness of solutions of WHVIs (or even WHGVIs)? Below, we answer the aforementioned question with an exceptional regularity condition expressed by \eqref{semi-ER} and some additional conditions.

\begin{theorem}\label{ER-1}
Given a nonempty closed subset $K$ in $C$ with $(K^\infty)^*$ being pointed, and two weakly homogeneous mappings {$f: C\rightarrow H$ and $g : H\rightarrow H$} defined by \eqref{weakhomo} with degrees $\delta_1> 0$ and $\delta_2> 0$, respectively. Suppose that $g^{-1}(C)\subseteq C$ and $\bar{g}({x})+{q}\in C$ for any $x\in H$ satisfying $\|x\|=1$, and one of the following conditions holds:
\begin{description}
  \item[\rm(i)] there exists no $(x,t)\in (H\setminus\{0\})\times\mathbb{R}_+$ such that
    \begin{eqnarray}\label{semi-ER}
     g^\infty(x)\in K^\infty,\;\; f^\infty({x})+t{x}\in (K^\infty)^{*}\;\; \text{\rm and}\;\;\langle f^\infty (x)+tx,g^\infty(x)\rangle=0,
    \end{eqnarray}
  and there exists a vector $d\in \mbox{\rm int}(K^\infty)$ such that
    \begin{eqnarray}\label{equ-10}
    \text{SOL}(x, g^\infty+d, K^\infty)=\{0\}=\text{SOL}(x, g^\infty, K^\infty);
    \end{eqnarray}
  \item[\rm(ii)] there exists no $(x,t)\in (H\setminus\{0\})\times\mathbb{R}_+$ such that \eqref{semi-ER} holds,
  and there exists a vector $d\in \mbox{\rm int}(K^\infty)$ such that
    \begin{eqnarray*}
    \text{SOL}(x, g+d, K^\infty)=\{0\}=\text{SOL}(x, g, K^\infty).
    \end{eqnarray*}
\end{description}
Then, WHGVI$(f,g,K)$ and WHGCP$(f,g,K^\infty)$ have nonempty compact solution sets.
\end{theorem}

\begin{proof}
First, we show that the result holds under condition (i). {Since (\ref{semi-ER}) implies that SOL$(f^{\infty},g^{\infty},K^\infty)=\{0\}$, we only need to prove {\rm ind}$\left((F^{\infty},g^{\infty})^{nat}_{K^\infty}(x),0\right)\neq 0$ where $F^\infty$ is any given continuous extensions of $f^\infty$ from Theorem \ref{Theorem-2.3}. For this purpose, we construct the following homotopy mapping:
$$
\mathcal{H}(x,t):=[(1-t)F^\infty({x})+tx]-\Pi_{(K^\infty)^*}
\{[(1-t)F^\infty(x)+tx]-
(g^\infty(x)+td)\}.
$$
Denote $\mathbb{Z}:=\{x\in H\mid \mathcal{H}(x,t)={0}\;\; \mbox{\rm for some}\ t\in [0,1]\}$. We show that $\mathbb{Z}$ is uniformly bounded. For the sake of contradiction, we assume that $\mathbb{Z}$ is not uniformly bounded. Then, we can find sequences $\{t_k\}\subseteq [0,1]$ and $\{{x^{k}}\}\subseteq H$ with $||{x^{k}}||\rightarrow \infty$ as ${k}\rightarrow \infty$ such that $\mathcal{H}({x^{k}},t_{{k}})=0$ for any ${k}$. It follows from $\mathcal{H}({x^{k}},t_{{k}})=0$ that
\begin{eqnarray*}
(1-t_{{k}})F^\infty(x^{k})+t_{k}x^k=\Pi_{(K^\infty)^*}
\{[(1-t_{k})F^\infty(x^k)+t_{k}x^k]-
[g^\infty (x^k)+t_{k}{d}]\},
\end{eqnarray*}
which, together with $(K^\infty)^*$ is a cone, implies that
\begin{eqnarray}\label{*}
\left\{\begin{array}{l}
(1-t_{k})F^\infty(x^k)+t_{k}x^k\in (K^\infty)^*,\quad g^\infty (x^k)+t_{k}{d}\in K^\infty,\\
\mbox{\rm and}\quad\langle g^\infty (x^k)+t_{k}{d}, [(1-t_{k})F^\infty(x^k)+t_{k}x^k]\rangle=0.
\end{array}\right.
\end{eqnarray}
Without loss of generality, we assume that $\delta_1>1$ (the proofs of case of $\delta_1\leq1$ can adopt similar procedure as the following), $\lim_{k\rightarrow\infty}t_{k}=\overline{t}$ and
$\lim_{k\rightarrow\infty}\frac{x^k}{\|x^k\|}=\overline{x}$ by the boundedness of sequences $\{t_{k}\}$ and $\{\frac{x^k}{\|x^k\|}\}$.

\begin{description}
\item[(I)] If $\overline{t}\neq1$, then $1-\overline{t}>0$. Dividing the equality in \eqref{*} by $\|x^k\|^{\delta_1+\delta_2}$, we obtain that
\begin{eqnarray}\label{2}
\left\langle\frac{g^\infty(x^k)+t_{k}{d}}
{\|x^k\|^{\delta_2}}, \frac{(1-t_{k})F^\infty(x^k)+t_{k}
x^k}{\|x^k\|^{\delta_1}}\right\rangle= 0.
\end{eqnarray}
Furthermore, since $\delta_1>1$, by letting $k\rightarrow\infty$ in \eqref{2}, we have that $\langle g^\infty(\overline{x}),
(1-\overline{t})F^\infty(\overline{x})\rangle=0$.
Noting that $\|\overline{x}\|=1$, then $\bar{g}(\overline{x})+{q}\in C$ and $g(\overline{x})=g^\infty(\overline{x})+\bar{g}(\overline{x})+{q}\in C$, which implies that $\overline{x}\in C$ by the condition $g^{-1}(C)\subseteq C$, and then, $F^\infty(\overline{x})=f^\infty(\overline{x})$. Thus, we can obtain
\begin{eqnarray}\label{SOL-2}
\langle g^\infty(\overline{x}),
(1-\overline{t})f^\infty(\overline{x})\rangle=0.
\end{eqnarray}
Moreover, it is easy to see that

\begin{eqnarray}\label{SOL-3}
\left\{\begin{array}{l}
g^\infty(\overline{x})=\lim_{k\rightarrow\infty}\frac{g^\infty(x^k)+t_kd}{\|x^k\|^{\delta_2}}\in K^\infty,\quad\mbox{\rm and}\\
(1-\overline{t})f^\infty(\overline{x})=\frac{(1-{t}_k)f^\infty(x^{k})+t_kx^k}{\|x^k\|^{\delta_1}}\in (K^\infty)^*.
\end{array}\right.
\end{eqnarray}
Let ${\overline{y}}=\sqrt[\delta_1]{1-\overline{t}}\cdot{\overline{x}}$, then ${\overline{y}}\neq {0}$; and it follows \eqref{SOL-2} and \eqref{SOL-3} that
$$
f^\infty (\overline{y})=(1-\overline{t})f^\infty(\overline{x})\in (K^\infty)^*,\quad g^\infty(\overline{y})={(1-\overline{t})}^{\frac{\delta_2}{\delta_1}}g^\infty(\overline{x})\in K^\infty,
$$
and
$$
\langle g^\infty(\overline{y}),F^\infty(\overline{y})\rangle=0.
$$
This is a contradiction to \eqref{semi-ER}!
\item[(II)] If $\overline{t}=1$, then dividing the equality of (\ref{*}) by $||x^k||^{\delta_2+1}$, we obtain that
$$
\left\langle\frac{g^\infty(x^k)+t_{k}{d}}
{\|x^k\|^{\delta_2}}, \frac{(1-t_{k})F^\infty(x^k)+t_{k}
x^k}{\|x^k\|}\right\rangle=0,$$
i.e.,
\begin{eqnarray}\label{1.1}
\left\langle\frac{g^\infty(x^k)+t_{k}d}
{\|x^k\|^{\delta_2}}, \frac{(1-t_{k})\|x^k\|^{\delta_1-1}F^\infty(x^k)}{\|x^k\|^{\delta_1}}+\frac{t_{k}
x^k}{\|x^k\|}\right\rangle=0.
\end{eqnarray}
Furthermore, we consider the following three cases.
\begin{itemize}
\item[$\bullet$] If $\lim_{k\rightarrow\infty}(1-t_{k})\|x^k\|^{\delta_1-1}=+\infty$, noting  $\lim_{k\rightarrow\infty}\frac{(1-t_{k})\|x^k\|^{\delta_1-1}}{\|x^k\|^{\delta_1-1}}=0$, then there exists some $N\in (0,{\delta_1-1})$ and $c>0$ such that
    $$
    \lim_{k\rightarrow\infty}\frac{(1-t_{k})\|x^k\|^{\delta_1-1}}{\|x^k\|^N}=c.
    $$
 In addition, it follows from \eqref{1.1} that
\begin{eqnarray*}
\left\langle\frac{g^\infty(x^k)+t_{k}d}
{||x^k||^{\delta_2}}, \frac{(1-t_{k})\|x^k\|^{\delta_1-1}F^\infty(x^k)}{\|x^k\|^{\delta_1}\cdot\|{x^k}\|^N}+\frac{t_{k}
x^k}{||x^k||^{N+1}}\right\rangle=0.
\end{eqnarray*}
Noting that $F^\infty(\overline{x})=f^\infty(\overline{x})$, then, letting $k\rightarrow\infty$, we have that
\begin{eqnarray}\label{$}
\langle g^\infty(\overline{x}),cf^\infty(\overline{x})\rangle=0,
\end{eqnarray}
and hence, follow the steps in (I), we can obtain that \eqref{$} is a contradiction to \eqref{semi-ER}!
\item[$\bullet$] If $\lim_{k\rightarrow\infty}(1-t_{k})\|x^k\|^{\delta_1-1}=0$, then letting $k\rightarrow\infty$ in \eqref{1.1}, we have that
\begin{eqnarray}\label{SOL-2-1}
\langle g^\infty(\overline{x}),
{\overline{x}}\rangle=0.
\end{eqnarray}
Moreover, it is easy to see that $g^\infty(\overline{x})\in K^\infty$ and
\begin{align*}
{\overline{x}}&=\lim_{k\rightarrow\infty}
[\frac{(1-t_{k})\|x^k\|^{\delta_1-1}F^\infty(x^k)}{\|x^k\|^{\delta_1}}+\frac{t_{k}x^k}
{\|x^k\|}]\\
&=\lim_{k\rightarrow\infty}
\frac{(1-t_{k})F^\infty(x^k)+t_{k}x^k}
{\|x^k\|}\in (K^\infty)^*,
\end{align*}
thus ${0}\neq {\overline{x}}\in\text{SOL}(x, g^\infty, K^\infty)$, which contradicts the second equality in \eqref{equ-10}!
\item[$\bullet$] If $\lim_{k\rightarrow\infty}(1-t_{k})\|x^k\|^{\delta_1-1}=c\;(c>0)$, then letting $k\rightarrow\infty$ in \eqref{1.1}, we have $\langle g^\infty(\overline{x}),cF^\infty(\overline{x})+{\overline{x}}\rangle=0$, which, together with $F^\infty(\overline{x})=f^\infty(\overline{x})$, implies that
\begin{eqnarray}\label{SOL-2-2}
\langle g^\infty(\overline{x}),
cf^\infty(\overline{x})+{\overline{x}}\rangle=0.
\end{eqnarray}
Moreover, we can easily obtain that $g^\infty(\overline{x})\in K^\infty$ and
\begin{eqnarray}\label{SOL-2-3}
cf^\infty(\overline{x})+\overline{x}=\lim_{k\rightarrow\infty}
\frac{(1-t_{k})F^\infty(x^k)+t_{k}x^k}
{\|x^k\|}\in (K^\infty)^*.
\end{eqnarray}
Since $c>0$, it follows from \eqref{SOL-2-2} and \eqref{SOL-2-3} that
$$\begin{array}{l}
g^\infty(\overline{x})\in K^\infty,\quad f^\infty(\overline{x})+\frac{1}{c}\overline{x}\in (K^\infty)^*\quad
\mbox{\rm and}\quad \langle g^\infty(\overline{x}),
f^\infty(\overline{x})+\frac{1}{c}{\overline{x}}\rangle=0,
\end{array}$$
which is a contradiction to \eqref{semi-ER}!
\end{itemize}
\end{description}
Thereby, by combining (I) with (II), we have that $\mathbb{Z}$ is uniformly bounded.}

Now, let $\Omega'$ be a bounded open set in $H$, which contains $\mathbb{Z}$, we have that $0\notin \mathcal{H}(\partial\Omega',t)$ for any $t\in [0,1]$. By the homotopy invariance principle of the degree, we have that,
\begin{eqnarray}\label{ind1=ind0}
\text{ind}\left((F^{\infty},g^{\infty})^{nat}_{K^\infty}(x),0\right)=\text{deg}(\mathcal{H}(x,0),
\Omega',0)=\text{deg}(\mathcal{H}(x,1),\Omega',0).
\end{eqnarray}

Noting that when $x$ is close to ${0}$, it holds that $x-[g^\infty(x)+d]$ is close to $-d$. In addition, since {$\Pi_{(K^\infty)^{*}}(-d)=0$ by $d\in \text{int}(K^\infty)$}, when $x$ is close to ${0}$, we have $\mathcal{H}(x,1)=x-{0}=x$,
which means that deg$(\mathcal{H}(x,1),\Omega^{'},{0})=1$ where $\Omega^{'}$ is a small open neighborhood of ${0}$. Therefore,
$$
\text{deg}(\mathcal{H}(x,1),\Omega,{0})=\text{ind}(\mathcal{H}(x,1),0)=\text{deg}(\mathcal{H}(x,1),\Omega^{'},{0})=1
$$
from the first equality in \eqref{equ-10}. Furthermore, from \eqref{ind1=ind0}, ind$\left((F^{\infty},g^{\infty})^{nat}_{K^\infty}(x),0\right)=1$.

Second, by constructing the homotopy mapping by
$$
\mathcal{H}(x,t):=[(1-t)F^\infty({x})+tx]-\Pi_{(K^\infty)^*}
\{[(1-t)F^\infty(x)+tx]-
(g(x)+td)\}
$$
and repeating the above procedure, we can obtain that ind$\left((F^{\infty},g^{\infty})^{nat}_{K^\infty}(x),0\right)=1$ when (ii) holds. This completes the proof.
\end{proof}

\begin{remark}
\rm{(i)} Suppose that $g(0)=0$, then the conditions described by the distance inequalities in $(a)$ of Theorems \ref{Theorem 2.4} and \ref{Theorem 2.5} hold naturally.

\rm{(ii)} If $f$ is strictly monotone with respective to $g$ on $H$, then the condition that $g$ is an injective mapping in $(a)$ of Theorem \ref{Theorem 2.5} holds naturally.
\end{remark}

\begin{remark}
\rm{(i)} When $g(x)=x$, it is obvious that $g^{-1}(C)\subseteq C$ and $\bar{g}({x})+{q}\in C$ for any $x\in H$ satisfying $\|x\|=1$, and the condition given in \eqref{equ-10} vanishes.

\rm{(ii)} In Theorem \ref{Theorem 2.4}, an essential condition is that $K^\infty$ is pointed; while Theorem \ref{ER-1} requires that $(K^\infty)^*$ is pointed. As there are a lot of cases that $K^\infty$ is pointed but $(K^\infty)^*$ is not, or $(K^\infty)^*$ is pointed but $K^\infty$ is not, Theorem \ref{Theorem 2.4} and Theorem \ref{ER-1} are two different results which are complements to each other.

\rm{(iii)} Here, we recall that a tensor ${\cal A}\in\mathbb{R}^{[m,n]}$ is said to have the $R$-property (see \cite{G17}) if SOL$({\cal A}x^{m-1},x,\mathbb{R}^n_+)=\{0\}=$SOL$({\cal A}x^{m-1}+d,x,\mathbb{R}^n_+)$ for some $d>0$; and it is called an $ER$-tensor (see \cite{WHB16}) if there exists no $(x,t)\in (\mathbb{R}^n_+\backslash\{0\})\times\mathbb{R}_+$ such that $({\cal A}x^{m-1})_i+tx_i=0$ if $x_i>0$, and $({\cal A}x^{m-1})_i\geq0$ if $x_i=0$. When the considered problems reduce to PCPs,  Theorem \ref{Theorem 2.4} reduces to Theorem 5.1 in \cite{G17}; while Theorem \ref{ER-1} reduces to Theorem 3.1 in \cite{LHL-18} with the involved leading tensors being $ER$-tensors.

\rm{(iv)} Let $d=(1,\cdots,1)^\top$ in the definition of the $R$-property, that is, if there exists no $(x,t)\in (\mathbb{R}^n_+\backslash\{0\})\times\mathbb{R}_+$ such that $({\cal A}x^{m-1})_i+t=0$ if $x_i>0$, and $({\cal A}x^{m-1})_i+t\geq0$ if $x_i=0$, then, $\cal{A}$ is called an $R$-tensor (see \cite{SQ17}). When the considered problems reduce to TCPs, Theorem \ref{Theorem 2.4} and Theorem \ref{ER-1} reduce to \cite[Theorem 4.1]{WHB16} and \cite[Theorem 4.2]{WHB16}, respectively.
\end{remark}

\begin{corollary}\label{cor-add4}
Suppose that all the conditions in Theorem \ref{ER-1} are satisfied. If additionally $f$ is strictly monotone with respect to $g$ on $K (K^\infty)$, then WHGVI$(f,g,K)$ (WHGCP$(f,g,K^\infty))$ has a unique solution.
\end{corollary}

Below, we construct an example satisfying all the conditions in Theorem \ref{ER-1}.
\begin{example}\label{exam-add0}
Let $H=C=\mathbb{R}^2$ and $K=\{(s,t)^\top\mid s\geq1,t\geq0\}$. Here, we consider WHGVI$(f,g,K)$, where $f:C\rightarrow H$ and $g:H\rightarrow H$ are defined as follows:
\begin{eqnarray*}
f(x)=\left(
\begin{array}{c}
x_1^{5}+x_1+1\\
x_2^{5}+x_2+1
\end{array}\right)\quad \mbox{\rm and}\quad
g(x)=\left(
\begin{array}{c}
x_1^{3}+x_1\\
x_2^{3}
\end{array}\right).
\end{eqnarray*}
\end{example}

Obviously, $g^{-1}(K)\subseteq C$ and $K^\infty=(K^\infty)^*=\mathbb{R}^2_+$. Below, we show that all the conditions in Theorem \ref{ER-1} are satisfied, and hence, WHGVI$(f,g,K)$ in Example \ref{exam-add0} has a unique solution.

\begin{itemize}
  \item It is not difficult to obtain that $g^{-1}(C)\subseteq C$ and $\bar{g}(x)+q\in C$ for any $x\in \mathbb{R}^2$.
  \item Since $g^\infty=(x_1^3,x_2^3)^\top$ and int$(K^\infty)=\mathbb{R}^2_{++}\neq\emptyset$, then, obviously, for any $d\in \mbox{\rm int}(K^\infty)$,
      $$
      \text{SOL}(x, g^\infty+d, K^\infty)=\{{0}\}{=\text{SOL}(x, g^\infty, K^\infty)}.
      $$
 \item It is easy to see that SOL$(f^\infty+tx, g^\infty, K^\infty)=\{{0}\}$ for any $t\geq0$.
  \item Besides, $f$ is strictly monotone with respect to $g$ on $K$.
\end{itemize}
Thus, by combining the above four cases, all the conditions in Theorem \ref{ER-1} hold, and hence, WHGVI$(f,g,K)$ in Example \ref{exam-add0} has a unique solution. In fact, it can be proved that $x^*=(x^*_1,0)^\top$ is the unique solution to WHGVI$(f,g,K)$, where $x^*_1$ satisfying $(x^*_1)^3+x^*_1-1=0$.

\section{Subcases of WHGVIs}\label{s4}
The WHGVI is a wide class of problems. When we restrict the set $K$ and/or mappings $g$ and $f$ to some special cases, we can obtain the corresponding uniquely solvable results of these problems. Actually, the results we obtain in Sections \ref{s3}, \ref{s6} and \ref{s9} can all reduce to the following subcases and get some nice results. In the following, we will not list all of these results, and only compare some of these results with the existing results.

\subsection{Reducing to special VIs}\label{s41}

In this subsection, we consider several VIs which are subclasses of WHGVIs, including WHVIs (i.e., weakly homogeneous variational inequalities), PGVIs (i.e., generalized polynomial variational inequalities), and PVIs (i.e., polynomial variational inequalities).

{\bf I. WHVIs}. Let $g(x)=x$, then, WHGVI$(f,g,K)$ reduces to the WHVI, denoted by WHVI$(f,K)$, which was studied by Gowda and Sossa \cite{GS19} and Ma et al. \cite{MZH19}. When reducing Theorem \ref{Theorem 2.4} from WHGVIs to WHVIs, we know that a relevant result is Theorem 5.1 in \cite{GS19}, which requires that $K^\infty$ is pointed.
When reducing Theorem \ref{ER-1} from WHGVIs to WHVIs, we can obtain the following result:
\begin{corollary}\label{ER-WHVI}
Given a nonempty closed subset $K$ in $C$ with $(K^\infty)^*$ being pointed and a weakly homogeneous mapping $f: C\rightarrow H$ defined by \eqref{weakhomo} with degree $\delta> 0$. Suppose that there exists no $(x,t)\in (H\setminus\{0\})\times\mathbb{R}_+$ such that $$x\in K^\infty,\quad f^\infty({x})+t{x}\in (K^\infty)^{*}\quad \mbox{\rm and}\quad \langle f^\infty (x)+tx,x\rangle=0.$$
Then, WHVI$(f,g,K)$ and WHCP$(f,g,K^\infty)$ have nonempty compact solution sets.
\end{corollary}

In addition, when reducing Theorem \ref{GWHVI} from WHGVIs to WHVIs, we can easily obtain the following result:

\begin{corollary}\label{WHVI}
Let $K$ be a nonempty closed convex set in $C$. Suppose that $f$ is a weakly homogeneous mapping defined by \eqref{weakhomo}. If there exists some $\theta\in K$ such that $\langle f(x)-f(\theta),x-\theta\rangle\geq0$ holds for any $x\in K$ and $\langle f^{\infty}(x),x\rangle\neq0$ for any $x\in B\bigcap K^{\infty}$, then, WHVI$(f,K)$ has a nonempty compact solution set.
\end{corollary}

Though many results on the nonemptiness and compactness of solution sets of WHVIs have been obtained by \cite{GS19} and \cite{MZH19}, we claim that Corollaries \ref{ER-WHVI} and \ref{WHVI} are new. Next, we only compare Corollary \ref{WHVI} with the main result in \cite{MZH19}.

From \eqref{weakhomo} we know that a weakly homogeneous mapping $f$ can be expressed as the following:
$$f(x)=f^{\infty}(x)+\bar{f}(x)+p=f^{\infty}(x)+\bar{f}(x)-\bar{f}(0)+p+\bar{f}(0).$$
Denote $\tilde{f}(x):=f^{\infty}(x)+\bar{f}(x)-\bar{f}(0)$ and $\tilde{p}:=p+\bar{f}(0)$, then $\tilde{f}^{\infty}=f^{\infty}$ and $\tilde{f}(0)=0$. More recently, \cite{MZH19} investigated the nonemptiness and compactness of the solution set of WHVI$(\tilde{f},K,\tilde{p})$. Although the expressions of WHVI$(\tilde{f},K,\tilde{p})$ in \cite{MZH19} and WHVI$(f,K)$ in this paper are different, there is a one-to-one correspondence between their solutions. The following is the main result in \cite{MZH19}.

\begin{theorem}\label{MZHvi}{\rm (\cite{MZH19})}
Let $K$ be a nonempty closed convex set in $C$, $\tilde{f}:C\rightarrow H$ be a weakly homogeneous mapping with degree $\gamma$, and $\tilde{p}\in H$. Suppose that the following conditions hold:
\begin{description}
  \item[\rm (i)] $\tilde{f}$ is $\eta$-copositive on $K$, that is, there exists a vector $\eta\in H$ such that $\langle \tilde{f}(x)-\eta,x\rangle\geq0$ holds for all $x\in K$; and
  \item[\rm (ii)] there exists a vector $\hat{x}\in K$ such that $\langle f(x),\hat{x}\rangle\leq0$ for all $x\in K$; and
  \item[\rm (iii)] let $S:=\text{SOL}(\tilde{f}^{\infty},K^\infty,0)$, where SOL$(\tilde{f}^{\infty},K^\infty,0)$ denotes the solution set of the WHVI$(\tilde{f}^{\infty},K^\infty,0)$, and $\tilde{p}+\eta\in \text{int}(S^*)$.
\end{description}
Then, WHVI$(\tilde{f},K,\tilde{p})$ has a nonempty compact solution set.
\end{theorem}

In the following, we use two examples to illustrate that the conditions in Corollary \ref{WHVI} are different from the conditions in Theorem \ref{MZHvi}.
\begin{example}\label{Ma-no}
Let $H=C=\mathbb{R}^2$ and $K=\{(s,t)^\top\mid s\geq-1,t\geq0\}$. We consider WHVI$(f,K)$, where $f(x)=(x_1^3-x_1^2,x_2^3+x_2)^\top$ is a weakly homogeneous mapping with degree 3 from $\mathbb{R}^2$ to $\mathbb{R}^2$.
\end{example}

We show that, for Example \ref{Ma-no}, all the conditions in Corollary \ref{WHVI} are satisfied, but at least one of the conditions in Theorem \ref{MZHvi} is not satisfied.

From Example \ref{Ma-no}, it is easy to see that $\theta=(1,0)^\top\in K$ and for any $x\in K$,
$$
\langle f(x)-f(\theta),x-\theta\rangle=x_1^2(x_1-1)^2+x_2^2(x_2^2+1)\geq0.
$$
Besides, it is also obvious that $\langle f^{\infty}(x),x\rangle=x_1^4+x_2^4\neq0$ for any $x\neq0$. Thus, $f$ satisfies all the conditions in Corollary \ref{WHVI}. However, $f=\tilde{f}$ is not $\eta$-copositive on $K$. Suppose on the contrary that there exists a vector $\eta=(\eta_1,\eta_2)^\top\in \mathbb{R}^2$ such that $\langle f(x)-\eta,x\rangle\geq0$ holds for all $x\in K$. Then,
$$
\left(\begin{array}{c}
x_1^3-x_1^2-\eta_1 \\ x_2^3+x_2-\eta_2
\end{array}\right)
\left(\begin{array}{c}
x_1 \\ x_2
\end{array}\right)\geq0,\quad \forall x_1\geq-1,x_2\geq0.
$$
Let $x_2=0$, then, we have $(x_1^3-x_1^2-\eta_1)x_1\geq0$. Now we consider three cases:
\begin{itemize}
  \item if $\eta_1\geq 0$, for $x_1\in(0,1)$ we have $x_1^3-x_1^2-\eta_1<0$, thus $(x_1^3-x_1^2-\eta_1)x_1<0$;
  \item if $\eta_1\in(-1,0)$, let $x_1=\eta_1/10<0$, then we have
      $$
      x_1^3-x_1^2-\eta_1=\frac{\eta_1^3}{1000}-\frac{\eta_1^2}{100}-\eta_1=\frac{\eta_1}{1000}(\eta_1^2-10\eta_1-1000)>0,
      $$
      thus $(x_1^3-x_1^2-q_1)x_1<0$;
  \item if $\eta_1\leq -1$, for $x_1\in(-1/2,0)$ we have $x_1^3-x_1^2-\eta_1>0$, thus $(x_1^3-x_1^2-\eta_1)x_1<0$.
\end{itemize}
Therefore, $f$ is not $\eta$-copositive on $K$. This indicates that at least one of the conditions in Theorem \ref{MZHvi} is not satisfied.

\begin{example}\label{Ma-yes}
Let $H=\mathbb{R}^2$, $C=\mathbb{R}^2_+$ and $K=\{(s',t')^\top\mid s'\geq0,t'\in[0,2\pi]\}$. We consider WHVI$(f,K)$, where $f(x)=(x_1+\sin x_1+1,\sin x_2+2)^\top$ is a weakly homogeneous mapping with degree 1 from $\mathbb{R}^2_+$ to $\mathbb{R}^2$.
\end{example}

We show that, for Example \ref{Ma-yes}, all the conditions in Theorem \ref{MZHvi} are satisfied, but at least one of the conditions in Corollary \ref{WHVI} is not satisfied.

From Example \ref{Ma-yes} it is easy to see that $\tilde{f}(x)=(x_1+\sin x_1,\sin x_2)^\top$ and $\tilde{p}=p=(1,2)^\top$. Let $\eta=(-1,-1)^\top$, then for any $x\in K$ we have
$$
\langle \tilde{f}(x)-\eta,x\rangle=x_1^2+x_1(1+\sin x_1)+x_2(1+\sin x_2)\geq0.
$$
Thus, $\tilde{f}$ is $\eta$-copositive on $K$. Let $\hat{x}=(0,0)^\top\in K$, then, $\langle f(x),\hat{x}\rangle\leq0$ holds for any $x\in K$. Besides, since $\tilde{f}^{\infty}(x)=(x_1,0)^\top$, we have $S:=\text{SOL}(\tilde{f}^{\infty},K^{\infty},0)=\{(0,0)^\top\}$ and $S^*=\mathbb{R}^2$. Furthermore, we have $\tilde{p}+\eta=(1,2)^\top+(-1,-1)^\top=(0,1)^\top\in\text{int}(S^*)$. Thus, the three conditions in Theorem \ref{MZHvi} hold. However, we cannot find a vector $\theta\in K$ such that $\langle f(x)-f(\theta),x-\theta\rangle\geq0$ holds for any $x\in K$. Suppose on the contrary there exists such a vector $\theta=(\theta_1,\theta_2)^\top\in K$. Here, we consider two cases:
\begin{itemize}
  \item Suppose that $\theta_1\geq0$ and $\theta_2\in[0,\pi]$. Let $x_1=\theta_1$, then, there exists $x_2>\theta_2$ which satisfies $\sin x_2<\sin\theta_2$, and hence,
       $$
       \langle f(x)-f(\theta),x-\theta\rangle=(\sin x_2-\sin\theta_2)(x_2-\theta_2)<0;
       $$
  \item Suppose that $\theta_1\geq0$ and $\theta_2\in(\pi,2\pi]$. Let $x_1=\theta_1$, then, there exists $x_2<\theta_2$ which satisfies $\sin x_2>\sin\theta_2$, and hence,
       $$
       \langle f(x)-f(\theta),x-\theta\rangle=(\sin x_2-\sin\theta_2)(x_2-\theta_2)<0.
       $$
\end{itemize}
Thus, there does not exist a vector $\theta\in K$ such that $\langle f(x)-f(\theta),x-\theta\rangle\geq0$  holds for any $x\in K$. This means that, for Example \ref{Ma-yes}, at lest one of the conditions in Corollary \ref{WHVI} is not satisfied.

\begin{remark}
It should be noted that, in both \cite{GS19} and \cite{MZH19}, the authors do not investigate the uniqueness of solutions to WHVI$(f,K)$. However, the uniqueness results in Sections \ref{s3}, \ref{s6} and \ref{s9} can all reduce to WHVIs. Thus, our uniqueness results in effect enrich the diversity of the results in WHVIs.
\end{remark}

{\bf II. PGVIs}. Let $H=\mathbb{R}^n$, and $f$ and $g$ be two polynomials defined by \eqref{poly} from $\mathbb{R}^n$ to $\mathbb{R}^n$, then, the WHGVI reduces to the PGVIs studied by Wang et al. \cite{WHX19}, which is denoted by GPVI$(f,g,K)$.


Reducing from WHGVIs to GPVIs, by Theorem \ref{main}, we immediately obtain the following result.

\begin{corollary}\label{c4}
Let $K$ be a nonempty closed convex set in $\mathbb{R}^n$, $f,g:\mathbb{R}^n\rightarrow \mathbb{R}^n$ be two polynomials defined by \eqref{poly}, and $\Omega^{\hat{x}}_r:=\{x\in H\mid \|g(x)\|<r\}$ where $r>\|\Pi_K(\hat{x})\|$ for any given $\hat{x}\in H$. Suppose that deg$(g(\cdot),\Omega^{\hat{x}}_r,\Pi_K(\hat{x}))$ is defined and nonzero, and the following conditions hold:
\begin{description}
  \item[\rm(i)]
  $f$ is strictly monotone with respect to $g$ on $K$; and
  \item[\rm(ii)]
  $\langle f^{\infty}(x),g^{\infty}(x)\rangle\neq0$ for any $x\in B\bigcap R$.
\end{description}
Then, GPVI$(f,g,K)$ has a unique solution.
\end{corollary}

Corollary \ref{c4} is a corrected version of Theorem 2 in \cite{WHX19}, since a restricted condition for the mapping $g$ (such as the degree condition in Theorem \ref{main}) was unnoticed in Theorem 2 in \cite{WHX19}.

\begin{remark}
We can obtain several uniquely solvable results from those in Sections \ref{s3}, \ref{s6} and \ref{s9} when WHGVIs reduce to WHVIs.
\end{remark}

{\bf III. PVIs}. Let $H=\mathbb{R}^n$, $g(x)=x$ and $f$ be a polynomial defined by \eqref{poly}, then, WHGVI$(f,g,K)$ reduces to the PVI studied by Hieu \cite{H18}. We denote it by PVI$(f,K)$. From Theorems \ref{main}, \ref{Theorem 2.4}, \ref{ER-1}, we can obtain some results for PVIs, in which the result reducing from Theorem \ref{Theorem 2.4} is exactly Theorem 5.1 in \cite{GS19}.

\begin{corollary}\label{pvi-1}
Given a nonempty closed convex subset $K$ of $\mathbb{R}^n$ and a polynomial $f$ defined by \eqref{poly}. Suppose that $f$ is strictly monotone on $K$ and $\langle f^\infty(x),x\rangle\neq0$ for any $x\in B\bigcap K^\infty$. Then, PVI$(f,K)$ has a unique solution.
\end{corollary}

\begin{corollary}\label{pvi-3}
Given a nonempty closed convex subset $K$ of $\mathbb{R}^n$ with $(K^\infty)^*$ being pointed and a polynomial $f$ defined by \eqref{poly}. Suppose that there exists no $(x,t)\in (\mathbb{R}^n\setminus\{0\})\times\mathbb{R}_+$ such that
\begin{eqnarray*}
x\in K^\infty,\;\; f^\infty(x)+t{x}\in (K^\infty)^{*}\;\; \text{\rm and}\;\;\langle f^\infty(x)+tx,x\rangle=0.
\end{eqnarray*}
Then, PVI$(f,K)$ has a nonempty compact solution set. Furthermore, if additionally $f$ is strictly monotone on $K$, then, PVI$(f,K)$ has a unique solution.
\end{corollary}

\begin{remark}
In \cite{H18}, the author gives a result on the existence and uniqueness of solutions to PVI$(f,K)$ under the condition that $0\in K$ and some additional conditions. However, Corollaries \ref{pvi-1} and \ref{pvi-3} do not require such a condition. Thus, when reducing from WHGVIs to PVIs, these corollaries enrich the theoretical results of the existence and uniqueness to solutions to PVIs.
\end{remark}

When $f(x)={\cal A}x^{m-1}+q$, where ${\cal A}\in \mathbb{R}^{[m,n]}$ and $q\in\mathbb{R}^n$, PVI$(f,K)$ further becomes the tensor variational inequality, denoted by TVI$({\cal A},K,q)$, investigated in \cite{WHQ18}. For TVI$({\cal A},K,q)$ we have the following result from Theorem \ref{main}.
\begin{corollary}\label{c3}
Given a nonempty closed convex subset $K$ of $\mathbb{R}^n$ and ${\cal A}\in \mathbb{R}^{[m,n]}$. Suppose that ${\cal A}x^{m-1}$ is strictly monotone on $K$, and $\langle {\cal A}x^{m-1},x\rangle\neq0$ for any $x\in B\bigcap K^{\infty}$. Then, TVI$({\cal A},K,q)$ has a unique solution for any given $q\in\mathbb{R}^n$.
\end{corollary}

If we further require that $0\in K$, then, it is easy to see that under the assumption that ${\cal A}x^{m-1}$ is strictly monotone on $K$, $\langle {\cal A}x^{m-1},x\rangle\neq0$ for any $x\in B\bigcap K^{\infty}$ is of course true. It is worth noting that this result is exactly Theorem 4.3 in \cite{WHQ18}.

\subsection{Reducing to CPs}\label{s42}
It is well-known that the CP is a subcase of the VI. Thus, in this subsection, we consider some cases of CPs which are subclasses of WHGVIs, including WHCPs (i.e., weakly homogeneous complementarity problems), PGCPs (i.e., generalized polynomial complementarity problems), and PCPs (i.e., polynomial complementarity problems).

{\bf I. WHCPs}. Let $C$ be a cone in $\mathbb{R}^n$ and $g(x)=x$, then, the WHGVI reduces to the WHCP, denoted by WHCP$(f,C)$. It is well-known that there is a uniqueness result in \cite{K71} for the conic complementarity problem. In the case of WHCPs, such a result states {\it the uniqueness of solutions to WHCP$(f,C)$ with the involved mapping $f:C\rightarrow H$ being strongly $C$-copositive on $C$, i.e., there exists a scalar $k > 0$ such that, for all $x \in C$, we have $\langle x, f(x)-f(0)\rangle \geq k\|x\|^2$}. According to the similar analysis in Propositions \ref{smne1} and \ref{smne2}, obviously, many weakly homogeneous mappings do not satisfy the strongly $C$-copositive property. Therefore, the uniqueness result in \cite{K71} cannot be directly applied to WHCPs in many cases. Below, we present some uniqueness results by reducing the relevant results from WHGVIs to WHCPs.

First, when reducing from WHGVIs to WHGCPs, we can obtain the following results by using some results in Section \ref{s6}.

\begin{theorem}\label{whgcp}
Let $C$ be pointed, and $f:C\rightarrow H$ and $g:H\rightarrow H$ be two weakly homogeneous mappings defined by \eqref{weakhomo} with degrees $\delta_1>0$ and $\delta_2>0$, respectively. Assume that $g(0)=0$. If $g$ is an injective mapping with $g^{-1}(C)\subseteq C$ and SOL$(f^{\infty},g^{\infty},C)=\{0\}$, then the following results are equivalent:
\begin{description}
\item[\rm ($a$)] WHGCP$(f+\xi,g,C)$ has a unique solution for any $\xi\in H$;
\item[\rm ($b$)] WHGCP$(f+\xi,g,C)$ has at most one solution for any $\xi\in H$.
\end{description}
\end{theorem}
\begin{proof}
First, we show that the result $(a)$ holds if the result $(b)$ holds. It is obvious that $0$ is a solution to WHGCP$(f^\infty+\bar{f}-\bar{f}(0)+d,g,C)$ where $d\in \text{int}(C^*)$. Let
$$F(x)=f(x)+\xi=f^\infty(x)+\bar{f}(x)+p+\xi,$$
then
$$F^\infty(x)+\bar{F}(x)-\bar{F}(0)=f^\infty(x)+\bar{f}(x)-\bar{f}(0).$$
Thus, from $(b)$ and the arbitrary of $\xi$, $0$ is the unique solution to WHGCP$(F^\infty+\bar{F}-\bar{F}(0)+d,g,C)$. According to Theorem \ref{Theorem 2.5}, we can obtain the uniqueness of WHGCP$(f+\xi,g,C)$.

Second, it is apparent that the result $(b)$ holds if the result $(a)$ holds.
\end{proof}

\begin{remark}
\rm{(i)} When the WHGCP reduces to the WHCP, the conditions about $g$ in Theorem \ref{whgcp} naturally hold, and the resulting conclusion for the WHCP is exactly Theorem 7.1 in \cite{GS19}. Thus, Theorem \ref{whgcp} is an extension of the uniqueness result obtained in \cite[Theorem 7.1]{GS19}.

\rm{(ii)} It follows from Lemma \ref{onesol} that WHCP$(f,C)$ has at most one solution when $f$ is strictly monotone with respect to $g$ on $C$. This and Theorem \ref{whgcp} will lead an existence and uniqueness result of WHCP$(f,C)$.

\rm{(iii)} We can obtain several uniquely solvable results from those in Sections \ref{s3}, \ref{s6} and \ref{s9} when WHGVIs reduce to WHCPs.
\end{remark}

Besides, in \cite{MZH19}, the authors also studied the properties of solution sets of WHCPs and obtained a nonemptiness and compactness result in \cite[Corollary 6.1]{MZH19} by reducing \cite[Theorem 3.1]{MZH19} (see Theorem \ref{MZHvi} in this paper) from WHVIs to WHCPs. To obtain the nonemptiness and compactness of solution sets to WHCPs, we can also reduce the relevant result in Corollary \ref{WHVI} from WHVIs to WHCPs and get a nonemptiness and compactness result for WHCPs.
\begin{corollary}\label{WHcp-add}
Given a weakly homogeneous mapping $f:C\rightarrow H$ defined by \eqref{weakhomo}. If there exists some $\theta\in C$ such that $\langle f(x)-f(\theta),x-\theta\rangle\geq0$ holds for any $x\in C$ and $\langle f^{\infty}(x),x\rangle\neq0$ for any $x\in B\bigcap C$, then, WHCP$(f,C)$ has a nonempty compact solution set.
\end{corollary}

Similar to Examples \ref{Ma-no} and \ref{Ma-yes}, it is easy to construct  examples to show that the result we obtain here is a different result from the one given by \cite{MZH19}.

{\bf II. PGCPs}. Let $H=\mathbb{R}^n$, $C$ be a cone in $\mathbb{R}^n$, and $f,g$ be two polynomials defined by \eqref{poly}, then, the WHGVI reduces to the PGCP studied by Ling et al \cite{LLH-18}. We denote it by GPCP$(f,g,C)$.
From Theorem \ref{main} we can directly obtain the uniqueness result:
\begin{corollary}\label{gpcp}
Let $C$ be a nonempty closed convex cone in $\mathbb{R}^n$, $f,g:\mathbb{R}^n\rightarrow \mathbb{R}^n$ be two polynomials defined by \eqref{poly}, and $\Omega^{\hat{x}}_r:=\{x\in H\mid \|g(x)\|<r\}$ where $r>\|\Pi_C(\hat{x})\|$ for any given $\hat{x}\in H$. Suppose that deg$(g(\cdot),\Omega^{\hat{x}}_r,\Pi_C(\hat{x}))$ is defined and nonzero, and the following conditions hold:
\begin{description}
  \item[\rm(i)] $f$ is strictly monotone with respect to $g$ on $C$; and
  \item[\rm(ii)] $\langle f^{\infty}(x),g^{\infty}(x)\rangle\neq0$ for any $x\in B\bigcap R$.
\end{description}
Then, GPCP$(f,g,C)$ has a unique solution.
\end{corollary}

More recently, in \cite{ZHM-19}, the authors also considered the GPCP and discussed the uniqueness of solutions to such a class of problems. One of the main result about the uniqueness given by \cite[Theorem 4.8]{ZHM-19} needs to find a vector $d\in\text{int}(C)\bigcap\text{int}(C^*)$. The following example shows that all the conditions of Corollary \ref{gpcp} are true, however, we may not find such a vector $d$.

\begin{example}\label{gpcp11}
Let $H=\mathbb{R}^n$ and $C=\{(t,2t)^\top\mid t\in\mathbb{R}\}$ be a cone. Here, we consider GPCP$(f,g,C)$, where $f,g:\mathbb{R}^n\rightarrow \mathbb{R}^n$ are defined as follows: $f(x)=(x_1^{3}+1,x_2^{3}+1)^\top$ and $g(x)=(x_1+1,x_2)^\top$.
\end{example}

From Example \ref{gpcp11} it is easy to see that:
\begin{itemize}
  \item $g$ is a bijective mapping on $\mathbb{R}^2$, thus, the degree condition about $g$ in Corollary \ref{gpcp} holds;
  \item $\langle f^\infty,g^\infty\rangle=x_1^4+x_2^4\neq0$ whenever $x\neq0$;
  \item for any $x=(x_1,2x_1)^\top\in C$ and $y=(y_1,2y_1)^\top\in C$ with $x\neq y$, we have $$\langle f(x)-f(y),g(x)-g(y)\rangle=17(x_1^3-y_1^3)(x_1-y_1)>0.$$
\end{itemize}
Thus, all the conditions in Corollary \ref{gpcp} hold. However, $C^*=\{(0,0)^\top\}$ leads to the fact that int$(C^*)$ is empty. Then, at least one condition of Theorem 4.8 in \cite{ZHM-19} is not satisfied.\vspace{2mm}

It is noticeable that Theorem \ref{Theorem 2.4} reduces to Theorem 4.8 in \cite{ZHM-19} when the WHGCP reduces to the GPCP. Thus, Example \ref{gpcp11} can also be used to verify that the uniqueness result given in Sections \ref{s3} and \ref{s6} are different from each other even reducing to GPCPs.
\vspace{2mm}

{\bf III. PCPs}. Let $H=\mathbb{R}^n$, $C=\mathbb{R}^n_+$, $g(x)=x$, and $f$ be a polynomial defined by \eqref{poly}, then, WHGVIs reduce to PCPs studied by \cite{G17}. We denote it by PCP$(f)$. In \cite{G17}, the author obtained a series of good results about the nonemptiness and compactness of solution sets to PCPs. A uniqueness result was also obtained in \cite[Theorem 6.1]{G17}. When the WHCP reduces to the PCP, Theorem \ref{whgcp} reduces to Theorem 6.1 in \cite{G17}.

Besides, \cite{LHL-18} investigated the nonemptiness and compact of the solution set, the uniqueness of solutions, and the error bounds of PCPs with the help of the structured tensors, where they gave a uniqueness result under the assumption that the involved mapping $f$ is an $m$-uniform $P$-function (i.e., there exists a constant $c>0$ such that $\max_{1\leq i\leq n}[x_i-y_i][f_i(x)-f_i(y)]\geq c\|x-y\|^m$ holds for any $x,y\in \mathbb{R}^n_+$).

Notice that when $g(x)=x$ and $f$ is a polynomial in WHGVI$(f,g,K)$, the exceptional regularity condition expressed by \eqref{semi-ER} reduces to the condition that the leading tensor of $f$ is an $ER$-tensor.
When WHGVIs reduce to PCPs, from Theorem \ref{ER-1} we can obtain the following result.

\begin{corollary}\label{pcp}
Let $f:\mathbb{R}^n\rightarrow \mathbb{R}^n$ be a polynomial defined by \eqref{poly}. Suppose that ${\cal A}^{(1)}$ is an $ER$-tensor and $f$ has $P$-property on $\mathbb{R}^n_+$. Then PCP$(f)$ has a unique solution.
\end{corollary}

\begin{remark}
\rm{(i)} In Theorem 3.1 in \cite{LHL-18}, the authors proved that if ${\cal A}^{(1)}$ is an $ER$-tensor, then the solution set of PCP$(f)$ is nonempty. This, together with Lemma \ref{onesol2}, implies that the result in Corollary \ref{pcp} holds.

\rm{(ii)} Since each of the condition that $f$ has $P$-property and the condition that ${\cal A}^{(1)}$ in $f$ is an $ER$-tensor is weaker than the condition that $f$ is an $m$-uniform $P$-function, it is sure that the first two conditions together are not stronger than the third condition. Therefore, Lemma 3.4 in \cite{LHL-18} can be seen as a corollary of Corollary \ref{pcp}.
\end{remark}

If $f(x)={\cal A}x^{m-1}+q$ for all $x\in \mathbb{R}^n$, where ${\cal A}\in\mathbb{R}^{[m,n]}$ and $q\in\mathbb{R}^n$, then, PCP$(f)$ reduces to the TCP. Like the TVI is an important subclass of the PVI, the TCP is a vital subcase of the PCP, which has attracted wide attention in recent years, and many papers consider the uniqueness of solutions to TCPs. When WHGVIs reduce to TCPs, Theorem \ref{Theorem 2.5} reduces to Theorem 3.7 in \cite{LLV-18}, which is derived based on Theorem 4.1 in \cite{BHW16}.

\section{Conclusion}\label{s5}

In this paper, with the help of the degree theory and the properties of weakly homogeneous mappings, we obtained several results on the unique solvability of WHGVIs, which were derived by making use of the exceptional family of elements for a pair of mappings, Karamardian-type theorems, and the exceptional regularity, respectively. In our main results, one of the main conditions is the strict monotonicity, which is weaker than the classical condition of strong monotonicity. Since the WHGVI provides a unified model for several classes of special VIs and CPs studied in recent years, this paper can be regarded as a unified treatment of the unique solvability of these subclasses in the sense that our conclusions can either reduce to known conclusions or give some new conclusions for these problems.

Up to now, the research on VIs and CPs with weakly homogeneous mappings mainly focuses on the nonemptiness and compactness of solution sets, and the unique solvability. One of the future issues is to study the theory of error bounds and the stability of solutions. Another future research topic is how to design efficient algorithms to solve these problems.

\end{document}